# MOMENTS AND TAILS IN MONOTONE-SEPARABLE STOCHASTIC NETWORKS

By François Baccelli[1] and Serguei Foss[2]

*INRIA-ENS, Institute of Mathematics and Heriot-Watt University*

A network belongs to the monotone separable class if its state variables are homogeneous and monotone functions of the epochs of the arrival process. This framework, which was first introduced to derive the stability region for stochastic networks with stationary and ergodic driving sequences, is revisited. It contains several classical queueing network models, including generalized Jackson networks, max-plus networks, polling systems, multiserver queues, and various classes of stochastic Petri nets. Our purpose is the analysis of the tails of the stationary state variables in the particular case of i.i.d. driving sequences. For this, we establish general comparison relationships between networks of this class and the $GI/GI/1/\infty$ queue. We first use this to show that two classical results of the asymptotic theory for $GI/GI/1/\infty$ queues can be directly extended to this framework. The first one concerns the existence of moments for the stationary state variables. We establish that for all $\alpha \geq 1$, the $(\alpha + 1)$-moment condition for service times is necessary and sufficient for the existence of the $\alpha$-moment for the stationary maximal dater (typically the time to empty the network when stopping further arrivals) in any network of this class. The second one is a direct extension of Veraverbeke's tail asymptotic for the stationary waiting times in the $GI/GI/1/\infty$ queue. We show that under subexponential assumptions for service times, the stationary maximal dater in any such network has tail asymptotics which can be bounded from below and from above by a multiple of the integrated tails of service times. In general, the upper and the lower bounds do not coincide. Nevertheless, exact asymptotics can be obtained along the same lines for various special cases of networks, providing direct extensions of Veraverbeke's tail asymptotic for the stationary waiting times in the $GI/GI/1/\infty$ queue. We

Received June 2001; revised November 2002.
[1]Supported in part by the TMR Alapedes project and INTAS.
[2]Supported in part by INTAS and Lyapunov Center.
*AMS 2000 subject classifications.* 60K25, 90B15, 60F10.
*Key words and phrases.* Queueing network, generalized Jackson network, ergodicity, subexponential random variable, tail asymptotics, Veraverbeke's theorem.

exemplify this on tandem queues (maximal daters and delays in stations) as well as on multiserver queues.

**1. Introduction.** We show in the present paper that properties which have been known for a long time for the tail asymptotics of isolated single server queues can be extended to the class of stochastic networks which are *monotone and separable*. This class, which was introduced in [6], contains several classical queueing network models like generalized Jackson networks, max-plus networks, polling systems and multiserver queues. This is also related to the class of topical (monotone and nonexpansive) maps of [18].

Section 2 summarizes the definition and main results that are known on this class of networks, and in particular the ergodic theorems that allow one to determine their stability region. The notion of maximal dater is recalled. In a generalized Jackson network, the maximal dater is the time to empty the network when stopping further arrivals. In a $G/G/1$ queue, this is just workload. In a FIFO tandem queue, this is end-to-end delay.

Section 3 focuses on the proof of the moment theorem. The assumptions that are needed here are limited to independence. We establish the following generalization of the classical $GI/GI/1$ queue moment theorem, which seems to be new within this setting: for all $\alpha \geq 1$, the $(\alpha+1)$-moment condition for service times in any monotone and separable network is necessary and sufficient for the existence of the $\alpha$-moment for the stationary maximal dater.

The subexponential tail asymptotic theorems are given in Sections 4 and 5. For surveys on the state of the art for this kind of asymptotics, see [20].

Section 4 gives generic upper and lower bounds which hold for all subexponential monotone separable networks and which only differ in the multiplicative constants.

Section 5 elaborates on the bounds established in Section 4. A corollary of Veraverbeke's theorem already proved in, for example, [1] and [2] states that, in the $GI/GI/1$ queue, large workloads occur on a typical event where a single large service time has taken place in the distant past, and all other service time are close to their mean. The main new result within our setting is Theorem 8 which extends the notion of typical event to subexponential monotone separable networks; large maximal daters occur when a single large service time has taken place in one of the stations and all other service time are close to their mean.

To the best of our knowledge, among the various classes of networks listed above, exact asymptotics are only known for irreducible max-plus networks [10]. The aim of Section 6 is to illustrate how the typical event theorem can be exploited to solve open questions on the exact asymptotics of other monotone separable networks. This is done for tandem queues in Section 6.1.1 and for multiserver queues in Section 6.2.

A first natural question is whether such asymptotics can be obtained for the maximal daters of all subexponential monotone separable networks. We have no general answer to this question yet. However, the choice of the tandem queue example to illustrate the potential use of the method was made on purpose; a tandem queue is both a generalized Jackson network and a reducible max-plus network. Exact asymptotics can be obtained along the same lines for the maximal daters of generalized Jackson networks and of reducible max-plus networks. These exact asymptotics require a lot of extra technical work, which go beyond the scope of the present paper and will be the object of two companion papers [7] and [8]. The exact asymptotics in polling systems is under investigation, too.

A second interesting question is whether such asymptotics can be extended to other characteristics than maximal daters. As it was shown in, for example, [21] maximal daters and individual waiting times may have fundamentally different asymptotics. This question is addressed in Section 6.1.2 where we show how to use the typical event theorem for monotone separable networks in order to derive the exact asymptotics for the stationary waiting or response times in individual queues of the tandem queue example.

## 2. Basic results on the monotone-separable networks.

2.1. *Framework.* Consider a stochastic network described by the following framework.

1. The network has a single input point process $N$, with points $\{T_n\}$; for all $m \leq n \in \mathbb{N}$, let $N_{[m,n]}$ be the $[m,n]$ restriction of $N$, namely the point process with points $\{T_l\}_{m \leq l \leq n}$.
2. The network has a.s. finite activity for all finite restrictions of $N$; for all $m \leq n \in \mathbb{N}$, let $X_{[m,n]}(N)$ be the time of the last activity in the network, when this one starts empty and is fed by $N_{[m,n]}$. We assume that for all finite $m$ and $n$ as above, $X_{[m,n]}$ is finite.

We assume that there exists a set of functions $\{f_l\}$, $f_l : \mathbb{R}^l \times K^l \to \mathbb{R}$, such that

$$(1) \qquad X_{[m,n]}(N) = f_{n-m+1}\{(T_l, \zeta_l), m \leq l \leq n\},$$

for all $n, m$ and $N$, where the sequence $\{\zeta_n\}$ is that describing service times and routing decisions.

We say that a network described as above is monotone separable if the functions $f_n$ are such that the following properties hold for all $N$.

1 (*Causality*). For all $m \leq n$,

$$X_{[m,n]}(N) \geq T_n.$$

2 (*External monotonicity*). For all $m \leq n$,
$$X_{[m,n]}(N') \geq X_{[m,n]}(N),$$
whenever $N' \stackrel{\text{def}}{=} \{T'_n\}$ is such that $T'_n \geq T_n$ for all $n$, a property which we will write $N' \geq N$ for short.

3 (*Homogeneity*). For all $c \in \mathbb{R}$ and for all $m \leq n$,
$$X_{[m,n]}(c+N) = X_{[m,n]}(N) + c.$$

4 (*Separability*). For all $m \leq l < n$, if $X_{[m,l]}(N) \leq T_{l+1}$, then
$$X_{[m,n]}(N) = X_{[l+1,n]}(N).$$

REMARK 1. Single-server queues, tandem queues, and generalized Jackson networks satisfy properties 1–4 above (see [5] and [6] for details).

REMARK 2. Using the terminology of the literature on idempotency (see, e.g., [18]), the monotone-separable framework can be rephrased in terms of so-called topical forms. Indeed, for all $m \leq n$, $X_{[m,n]}$ can be seen as function of the bi-infinite vector $(\ldots, T_{-2}, T_{-1}, T_0, T_1, T_2, \ldots)$ of $\mathbb{R}^{\mathbb{Z}}$. Since $X_{[m,n]} : \mathbb{R}^{\mathbb{Z}} \to \mathbb{R}$ is monotone and homogeneous, according to this terminology, the family $X_{[m,n]}$, $-\infty \leq m \leq n < \infty$, is a family of topical forms on $\mathbb{R}^{\mathbb{Z}}$. The link between these forms is established via the separability assumption, which allows one to study the asymptotic forms $X_{(-\infty,n]}$, which are the main objects of interest. Of particular interest to us are the statistical properties (moments, tail behavior, etc.) or the projective properties of the sequence $(T_n, X_{(-\infty,n]}) \in \mathbb{R}^2$.

2.2. *Maximal daters*. By definition, the $[m,n]$ maximal dater is
$$Z_{[m,n]}(N) \stackrel{\text{def}}{=} X_{[m,n]}(N) - T_n = X_{[m,n]}(N - T_n).$$
Note that $Z_{[m,n]}(N)$ is a function of $\{\zeta_l\}_{m \leq l \leq n}$ and $\{\tau_l\}_{m \leq l \leq n-1}$ only, where $\tau_n = T_{n+1} - T_n$. In particular, $Z_n(N) \stackrel{\text{def}}{=} Z_{[n,n]}(N)$ is not a function of $\{\tau_l\}_{-\infty < l < \infty}$.

LEMMA 1 (Internal monotonicity of $X$ and $Z$). *Under the above conditions, the variables $X_{[m,n]}$ and $Z_{[m,n]}$ satisfy the internal monotonicity property; for all $N$,*
$$X_{[m-1,n]}(N) \geq X_{[m,n]}(N),$$
$$Z_{[m-1,n]}(N) \geq Z_{[m,n]}(N), \qquad m \leq n.$$

In particular, the sequence $\{Z_{[-n,0]}(N)\}$ is nondecreasing in $n$. Put
$$Z \equiv Z_{(-\infty,0]} = \lim_{n \to \infty} Z_{[-n,0]}(N) \leq \infty.$$

LEMMA 2 (Subadditive property of $Z$). *Under the above conditions, $\{Z_{[m,n]}\}$ satisfies the following subadditive property: for all $m \leq l < n$, for all $N$,*

$$Z_{[m,n]}(N) \leq Z_{[m,l]}(N) + Z_{[l+1,n]}(N).$$

2.3. *Stochastic assumptions and main stability results.* Assume the variables $\{\tau_n, \zeta_n\}$ are random variables defined on a common probability space $(\Omega, \mathcal{F}, \mathbf{P}, \theta)$, where $\theta$ is an ergodic, measure-preserving shift transformation, such that $(\tau_n, \xi_n) \circ \theta = (\tau_{n+1}, \xi_{n+1})$. The following integrability assumptions are also assumed to hold:

$$\mathbf{E}\tau_n \stackrel{\text{def}}{=} \lambda^{-1} \stackrel{\text{def}}{=} a < \infty, \qquad \mathbf{E}Z_n < \infty.$$

We summarize the main results of [6].

LEMMA 3 (0–1 law). *Under the foregoing ergodic assumptions, either $Z = \infty$ a.s. or $Z < \infty$ a.s.*

The network is stable if $Z < \infty$ a.s. and unstable otherwise.

Denote by $Q = \{T'_n\}$ the degenerate input process with $T'_n = 0$ a.s. for all $n$.

LEMMA 4. *Under the foregoing ergodic assumptions, there exists a non-negative constant $\gamma(0)$ such that*

$$\lim_{n \to \infty} \frac{Z_{[-n,-1]}(Q)}{n} = \lim_{n \to \infty} \frac{\mathbf{E}Z_{[-n,-1]}(Q)}{n} = \gamma(0) \qquad a.s.$$

The main result on the stability region is the following theorem.

THEOREM 1. *If $\lambda\gamma(0) < 1$, then $Z < \infty$ a.s. If $Z < \infty$ a.s., then $\lambda\gamma(0) \leq 1$.*

2.4. *Further assumptions.* Most of the new results of the present paper will be obtained under the following independence assumption.

(IA). The sequences $\{\zeta_n\}$ and $\{\tau_n\}$ are mutually independent and each of them consists of i.i.d. random variables.

For certain results, we shall make the following additional assumption.

(AA). For all $i$,

(2) $$Z_i = Z_{[i,i]} = Y_i^{(1)} + \cdots + Y_i^{(r)},$$

where the r.v.'s $Y_i^{(j)}$ are nonnegative, independent of interarrival times, and such that the sequence of random vectors $(Y_i^{(1)}, \ldots, Y_i^{(r)})$ is i.i.d; general dependences between the components of the vector $(Y_i^{(1)}, \ldots, Y_i^{(r)})$ are allowed. In addition,

$$Z_{[n,0]}(Q) \geq \max_{j=1,\ldots,r} \sum_{i=n}^{0} Y_i^{(j)} \quad \text{a.s.} \tag{3}$$

2.5. *Upper and lower bound $G/G/1/\infty$ queues.* The results of this section are new. We assume stability, namely $\gamma(0) < a$. We pick an integer $L \geq 1$ such that

$$\mathbf{E} Z_{[-L,-1]}(Q) < La, \tag{4}$$

which is possible in view of Lemma 4. Without loss of generality, one can assume $T_0 = 0$.

To the input process $N$, we associate the following lower and upper bound processes: $N^- = \{T_n^-\}$, where, for all $k$ and $n$ in $\mathbb{Z}$ such that $n = (k-1)L + 1, \ldots, kL$, $T_n^- = T_{(k-1)L}$, and similarly, $N^+ = \{T_n^+\}$, where $T_n^+ = T_{kL}$ if $n = (k-1)L + 1, \ldots, kL$. Then for all $n$,

$$X_{[-n,0]}(N^-) \leq X_{[-n,0]}(N) \leq X_{[-n,0]}(N^+) \equiv Z_{[-n,0]}(N^+). \tag{5}$$

In other words, both upper and lower bound processes have batch arrivals (of size $L$).

Note that if (IA) holds, the r.v.'s $Z_{[-n,0]}(N^-) = X_{[-n,0]}(N^-) - T_{-L}$ and $Z_{[-n,0]}(N^+)$ have the same distribution and that the r.v.'s $Z_{[-n,0]}(N^-)$ and $T_{-L}$ are independent.

2.5.1. *Upper bound queue.* The next lemma, which establishes a first connection between monotone-separable networks and the $G/G/1/\infty$ queue, directly follows from the monotonicity and the separability assumptions.

LEMMA 5. *Assume $T_0 = 0$. For any $m < n \leq 0$,*

$$Z_{[m,0]}(N) \leq Z_{[n,0]}(N) + \max(0, Z_{[m,n-1]}(N) - \tau_{n-1}).$$

Put $Z_n = Z_{[n,n]}(N)$. Then the sequence $\{Z_n\}$ does not depend on $N$ and forms a stationary and ergodic sequence.

COROLLARY 1. *Assume $T_0 = 0$. For any $m < 0$,*

$$Z_{[m,0]} \equiv Z_{[m,0]}(N) \leq \max_{m \leq k \leq 0} \left( \sum_{i=k}^{0} Z_i - \sum_{i=k+1}^{0} \tau_i \right)$$

*with the convention $\sum_{1}^{0} = 0$.*

The main weakness of this upper bound comes from the fact that the corresponding queue may be unstable whereas the initial network is stable. This is taken care of by the upper bound described below.

COROLLARY 2. *The stationary maximal dater $Z \equiv Z_{(-\infty,0]}(N)$ is bounded from above by the stationary response time $\widehat{R}$ in the $G/G/1/\infty$ queue with service times*

$$\widehat{s}_n = Z_{[L(n-1)+1, Ln]}(Q) \tag{6}$$

*and interarrival times $\widehat{\tau}_n = T_{Ln} - T_{L(n-1)}$, where $L$ is the integer defined in (4). Since $\widehat{b} = \mathbf{E}\,\widehat{s}_n < \mathbf{E}\,\widehat{\tau}_n = La$, this queue is stable.*

PROOF. We have

$$Z = \lim_{n\to\infty} Z_{[-n,0]} = \lim_{k\to\infty} Z_{[-kL+1,0]} = \sup_{k\geq 0} Z_{[-kL+1,0]}$$

$$\leq \sup_{k\geq 0} Z_{[-kL+1,0]}(N^+) \leq \sup_{k\geq 0} \max_{-k\leq i\leq 0}\left(\widehat{s}_0 + \sum_{j=i}^{-1}(\widehat{s}_j - \widehat{\tau}_{j+1})\right)$$

$$= \widehat{s}_0 + \sup_{k\geq 0} \sum_{i=-k}^{-1}(\widehat{s}_i - \widehat{\tau}_{i+1}) = \widehat{R}.$$

In these relations, (5) was used to derive the first inequality, Corollary 1 was used in the last inequality; we also used the fact that

$$Z_{[L(n-1)+1, Ln]}(N^+) = Z_{[L(n-1)+1, Ln]}(Q)$$

and the convention $\sum_0^{-1} = 0$. □

The queue of Corollary 2 will be referred to as the *$L$-upper-bound $G/G/1/\infty$ queue* associated with the network.

Note that when (IA) holds, this queue is a $GI/GI/1/\infty$ queue. In this case, $\widehat{R} = \widehat{W} + \widehat{s}_0$, where $\widehat{W}$ is a stationary waiting time and $\widehat{W}$ and $\widehat{s}_0$ are independent.

Notice that under (AA),

$$\max_{j=1,\ldots,r} \sum_{i=L(n-1)+1}^{Ln} Y_i^{(j)} \leq \widehat{s}_n \leq \sum_{j=1}^{r} \sum_{i=L(n-1)+1}^{Ln} Y_i^{(j)} \quad \text{a.s.} \tag{7}$$

where the second inequality follows from the subadditive property of $Z$.

The following result does not require (AA) and holds for all monotone separable networks such that the sequence $\{Z_i\}$ is i.i.d.

We say that a nonnegative r.v. $X$ is *light tailed* if there exists a positive number $c$ such that $\mathbf{E} \exp(cX)$ is finite.

COROLLARY 3. *If $Z_0$ is light tailed and $\lambda\gamma(0) < 1$, then $Z_{(-\infty,0]}$ is light tailed too.*

PROOF. From Corollary 2, it is enough to prove that the response time $\widehat{R}$ in the $L$-upper-bound queue is light tailed. From well-known results on the $G/GI/1/\infty$ queue, in the stable case, the stationary response times $\widehat{R}$ are light tailed when the service times are light tailed. But from the subadditive inequality, we have

$$\widehat{s}_1 = Z_{[1,L]}(Q) \leq \sum_{i=1}^{L} Z_i, \tag{8}$$

which proves that $\widehat{s}_1$ is light tailed if $Z_0$ is. $\square$

2.5.2. *Lower bound fork-join queue.* The following result is immediate.

LEMMA 6. *Under Condition* (AA),

$$Z_{(-\infty,0]} \geq \underline{R} = \max_{j=1,\ldots,r} \sup_{n \leq 0} \left( \sum_n^0 Y_i^{(j)} - \sum_n^{-1} \tau_i \right). \tag{9}$$

The queue with service times $\{Y_i^{(j)}\}$ and interarrival times $\{\tau_i\}$ will be referred to as the $j$-lower-bound $G/G/1/\infty$ queue associated with the network. Let $R^{(j)}$ denote the stationary response time in this queue:

$$R^{(j)} = \sup_{n \leq 0} \left( \sum_n^0 Y_i^{(j)} - \sum_n^{-1} \tau_i \right).$$

Then the lower bound $\underline{R}$ defined in (9) is the stationary response time in the $r$-dimensional *fork-join queue* with service times $\{Y_i^{(j)}\}$, $j = 1,\ldots,r$ and interarrival times $\{\tau_i\}$.

2.6. *Examples.*

2.6.1. *Tandem queues.* Consider a stable $G/G/1/\infty \to \cdot/G/1/\infty$ tandem queue. Denote by $\{\sigma_n^{(i)}\}$ the sequence of service times in station $i = 1, 2$ and $\{\tau_n\}$ the sequence of interarrival times at the first station. Put $b^{(i)} = \mathbf{E}\sigma^{(i)}$, $a = \mathbf{E}\tau$ and $\rho^{(i)} = b^{(i)}/a < 1$. We have $\gamma(0) = \max(b^{(1)}, b^{(2)})$.

Tandem queues fall in the class of open Jackson networks, and in the class of open max-plus systems which both belong to the class of monotone separable networks (see below). We have the following representation for the maximal dater (see, e.g., [10])

$$
(10) \quad Z_{[-n,0]} = \sup_{-n \leq p \leq 0} \sup_{p \leq q \leq 0} \left( \sum_{m=p}^{q} \sigma_m^{(1)} + \sum_{m=q}^{0} \sigma_m^{(2)} - (T_0 - T_p) \right),
$$

$$
(11) \quad Z = Z_{(-\infty,0]} = \sup_{p \leq 0} \sup_{p \leq q \leq 0} \left( \sum_{m=p}^{q} \sigma_m^{(1)} + \sum_{m=q}^{0} \sigma_m^{(2)} - (T_0 - T_p) \right).
$$

Assumption (IA) is satisfied if the sequences $\{\tau_n\}$ and $\{\zeta_n \equiv (\sigma_n^{(1)}, \sigma_n^{(2)})\}$ are i.i.d. and mutually independent (we may allow a dependence between $\sigma_n^{(1)}$ and $\sigma_n^{(2)}$). Assumption (AA) is also satisfied here with $r = 2$ and $Y_n^{(i)} = \sigma_n^{(i)}$, $i = 1, 2$.

The maximal dater with index $n$ is the sojourn time of customer $n$ in the network, namely the time which elapses between its arrival in station 1 and its departure from station 2.

As for the $L$-upper-bound queue associated with this network, the expression for $\widehat{s}_n$ is here

$$
(12) \quad \widehat{s}_n = \max_{1 \leq j \leq L} \left( \sum_{i=1}^{j} \sigma_{(n-1)L+i}^{(1)} + \sum_{i=j}^{L} \sigma_{(n-1)L+i}^{(2)} \right).
$$

2.6.2. *Multiserver queues.* Let

$$W_n = (W_n^{(1)}, \ldots, W_n^{(m)})$$

be the Kiefer–Wolfowitz workload vector in the $GI/GI/m/\infty$ queue with interarrival times $\{\tau_n\}$ and service times $\{\sigma_n\}$. Here $n$ is the customer index and $W_n^{(i)}$, $i = 1, \ldots, m$, are the workloads of the servers at the $n$th arrival time, arranged in nondecreasing order. More precisely, we assume $W_0 = (0, \ldots, 0)$ and

$$
(13) \quad W_{n+1} = \mathbf{R}(W_n + \mathbf{e}_1 \sigma_n - \mathbf{i}\tau_n)^+
$$

for $i \geq 0$, where $\mathbf{e}_1 = (1, 0, \ldots, 0)$ and $\mathbf{i} = (1, 1, \ldots, 1)$ are $m$-dimensional vectors and the operator $\mathbf{R}$ permutes the components of a vector in nondecreasing order. For a multiserver queue, $\gamma(0) = \mathbf{E}\sigma_0/m$.

Assumption (IA) is satisfied under the assumption that the service times are i.i.d. Assumption (AA) is not satisfied here.

The maximal dater associated with customer $n$ is the time which elapses between its arrival and the time when all customers still present at its arrival time have left the system (including customer $n$),

$$Z_{[0,n]} = \max(W_n^{(1)} + \sigma_n, W_n^{(m)}).$$

2.6.3. *Generalized Jackson networks.* Consider a generalized Jackson network with $r$ stations. We denote by:

1. $\{\sigma_n^{(k)}\}$ the i.i.d. sequence of service times in station $k$.
2. $\{\mu_n^{(i)}\}$ the i.i.d. sequence of routing decisions from station $i$; with values in the set $\{1, \ldots, r\}$.
3. $\{\mu_n\}$ the i.i.d. sequence of routing decisions for the input process; also with values in the set $\{1, \ldots, r, r+1\}$, where $\mu_n^{(j)} = r+1$ means that a customer taking $n$th service at station $i$ leaves then the network.
4. $\{\tau_n\}$ the i.i.d. sequence of interarrival time.

Under these assumptions, both (IA) and (AA) are satisfied. We have

$$Y_1^{(j)} = \sum_1^{\nu(j)} \sigma_n^{(j)} \tag{14}$$

with $\nu(j)$ the total number of visits of customer 1 (the customer arriving at time $T_1$) to station $j$ in the $[1,1]$ restriction of the network, namely when this customer is the only one to enter the network. The random variables $\nu(j)$, $j = 1, \ldots, r$, are obtained from the sequences of routing decisions (see [5]).

In this case $Z_{[-n,0]}$ is the time which elapses between the arrival of customer 0 and the time when all customers have left the system, given that arrivals are stopped after $T_0$.

2.6.4. *Max-plus networks.* The class of open max-plus networks also falls in this framework (see, e.g., [6]). A typical example of this class is that of tandem queues. Tandem queues form a reducible open max-plus network. For examples of irreducible networks of this class, see [10].

**3. Integrability of stationary maximal daters.** We assume (IA) and stability, namely $\lambda \gamma(0) < 1$.

Let $\widehat{W}$ denote the stationary waiting time in the $L$-upper-bound $GI/GI/1/\infty$ queue of the network. The following result is well known.

LEMMA 7. *For any $\alpha > 1$, $\mathbf{E}\widehat{W}^{\alpha-1}$ is finite if and only if $\mathbf{E}\widehat{s}_0^\alpha$ is finite.*

Therefore, $\widehat{R} = \widehat{W} + \widehat{s}$ is such that $\mathbf{E}\widehat{R}^{\alpha-1}$ is finite if and only if $\mathbf{E}\widehat{s}_0^\alpha$ is finite.

COROLLARY 4. *If $\mathbf{E}Z_0^\alpha < \infty$, then $\mathbf{E}Z_{(-\infty,0]}^{\alpha-1} < \infty$.*

PROOF. We have

$$\widehat{s}_0 \leq \sum_{-L+1}^0 Z_i.$$

Therefore if $\mathbf{E}Z_0^\alpha < \infty$, then $\mathbf{E}\widehat{s}_0^\alpha$ is finite. Thus, $\mathbf{E}\widehat{W}^{\alpha-1}$ and $\mathbf{E}\widehat{R}^{\alpha-1}$ are finite, too. We conclude the proof by using the bound $Z_{(-\infty,0]} \leq \widehat{R}$ (see the proof of Corollary 2). $\square$

Under condition (AA), $\mathbf{E}Z_0^\alpha$ is finite if and only if for all $j$, $\mathbf{E}[(Z_0(j))^\alpha]$ is finite. The following theorem is then an immediate consequence of Lemmas 6 and 7.

THEOREM 2. *Under assumptions* (IA) *and* (AA), *if* $\mathbf{E}[Z_{(-\infty,0]}^{\alpha-1}]$ *is finite, so is* $\mathbf{E}Z_0^\alpha$.

EXAMPLES. All results are given under the assumption that the system under consideration is stable.

1. Tandem queues. The system response time has a moment of order $\alpha - 1$ iff the service times in both stations admit a moment of order $\alpha$.
2. Multiserver queues. In steady state, the time to empty the system has a moment of order $\alpha - 1$ if the service times admit a moment of order $\alpha$.
3. Generalized Jackson networks. The stationary maximal dater has a moment of order $\alpha - 1$ iff all service times have moments of order $\alpha$. Since the stationary maximal dater is not less than the residual workload at any station, we get if all service times have moments of order $\alpha$, then the stationary residual workloads and the stationary queue lengths at all stations have moments of order $\alpha - 1$. Since the number of customer services has a geometrical tail, one also deduces from this that stationary sojourn times also have moments of order $\alpha - 1$.

**4. Bounds for subexponential tail asymptotics.**

4.1. *Assumptions and notation.* Here and later in the paper, for strictly positive functions $f$ and $g$, the equivalence $f(x) \sim dg(x)$ with $d > 0$ means $f(x)/g(x) \to d$ as $x \to \infty$. This equivalence may also be rewritten as $f(x) = dg(x)(1+o(1)) = dg(x) + o(g(x)) = dg(x) + o(f(x))$, where $o(1)$ is a function which tends to 0 as $x$ tends to $\infty$, and $o(g(x))$ is a function such that $o(g(x))/g(x) \to 0$ as $x \to \infty$. By convention, the equivalence $f(x) \sim dg(x)$ with $d = 0$ means $f(x) = o(g(x))$. We will also use the following notation:

1. $f(x) = \Theta(g(x))$ to mean $\limsup f(x)/g(x) < \infty$ and $\liminf f(x)/g(x) > 0$,
2. $f(x) = O(g(x))$ to mean $\limsup f(x)/g(x) < \infty$.

4.1.1. *Tails.* Let $\xi$ be a nonnegative r.v. with distribution function $F$ such that $\mathbf{P}(\xi > x) \equiv 1 - F(x) \equiv \overline{F}(x) > 0$ for all $x$. Let $\xi_1, \xi_2$ be independent copies of $\xi$.

DEFINITION 1. $\xi$ has a *heavy-tailed* distribution (HT), if, for any $c > 0$,
$$\mathbf{E}\exp(c\xi) \equiv \int_0^\infty \exp(cx)\,dF(x) = \infty.$$

DEFINITION 2. $\xi$ has a *long-tailed* distribution (LT), if, for any $y > 0$,
$$\overline{F}(x+y) \sim \overline{F}(x) \qquad \text{as } x \to \infty.$$

Any LT distribution is HT.

DEFINITION 3. $\xi$ has a *subexponential* distribution (SE), if
$$\mathbf{P}(\xi_1 + \xi_2 > x) \sim 2\overline{F}(x) \qquad \text{as } x \to \infty.$$

Any SE distribution is LT. For basic properties of subexponential distributions, see, for example, [13].

4.1.2. *Network assumptions.* Consider a distribution function $F$ on $\mathbb{R}^+$ such that the following hold:

(a) $F$ is subexponential, with finite first moment $M = \int_0^\infty \overline{F}(u)\,du$, where $\overline{F}(u) = 1 - F(u)$ denotes the tail of $F$.

(b) The integrated tail distribution $F^s$,
$$F^s(x) = 1 - \min\left\{1, \int_x^\infty \overline{F}(u)\,du\right\} \equiv 1 - \overline{F}^s(x),$$
is subexponential.

Here are a few properties satisfied by $F$ that will be needed later on and that follow from the fact that $F^s$ is long tailed.

When $x \to \infty$,

(15) $$\overline{F}(x) = o(\overline{F}^s(x)).$$

As a corollary, there exists a nondecreasing integer-valued function $N_x \to \infty$ and such that, for all finite real numbers $b$,

(16) $$\sum_{n=0}^{N_x} \overline{F}(x + nb) = o(\overline{F}^s(x)), \qquad x \to \infty.$$

In particular,

(17) $$N_x \overline{F}(x) = o(\overline{F}^s(x)).$$

Such a c.d.f. $F$ being given, consider a monotone separable network satisfying (IA) and (AA) and such that the following equivalence holds when $x$ tends to $\infty$:

(c) For all $j = 1, \ldots, r$,

$$\mathbf{P}(Y_1^{(j)} > x) \sim d^{(j)} \overline{F}(x)$$

with $\sum_j d^{(j)} \equiv d > 0$.

For a monotone separable network, the three assumptions (a)–(c) will be referred to as (SE). Under (SE), the following holds:

(18) $$\int_x^\infty \mathbf{P}(Y_1^{(j)} > y) \, dy \sim d^{(j)} \overline{F}^s(x) \qquad \text{as } x \to \infty.$$

ASSUMPTION (H). We also introduce the following assumption (H):

(19) $$\mathbf{P}\left(\sum_1^r Y_1^{(j)} > x\right) \sim \mathbf{P}\left(\max_{1 \leq j \leq r} Y_1^{(j)} > x\right)$$
$$\sim \sum_1^r \mathbf{P}(Y_1^{(j)} > x) \sim \sum_1^r d^{(j)} \overline{F}(x).$$

Note that the very last equivalence follows from (SE). Assumption (H) is, for instance, satisfied in the particular case when the random variables $Y_1^{(j)}$ are mutually independent; in Section A.2, we give sufficient conditions for (H) to hold that go beyond this particular case.

Take any $1 \leq i_1, i_2 \leq r$, $i_1 \neq i_2$. Since

$$\mathbf{P}\left(\max_j Y_1^{(j)} > x\right) \leq \sum_j \mathbf{P}(Y_1^{(j)} > x) - \mathbf{P}(Y_1^{(i_1)} > x, Y_1^{(i_2)} > x),$$

we deduce from (19) that

(20) $$\mathbf{P}(Y_1^{(i_1)} > x, Y_1^{(i_2)} > x) = o(\overline{F}(x)).$$

REMARK 3. In what follows, we will not need i.i.d. assumptions on the interarrival times $\{\tau_n\}$. As it follows from Theorem 14, the results we prove will hold also in the more general situation when these variables satisfy the following three conditions:

1. $\{\tau_n\}$ forms a stationary ergodic sequence with a finite positive mean $\mathbf{E}\tau_1 = a$.
2. $\{\tau_n\}$ is independent of $\{Y_n^{(j)}, j = 1, \ldots, r\}$.
3. For all $\widetilde{a} < a$,

$$\mathbf{P}\left(\sup_{n \geq 0}\left(n\widetilde{a} - \sum_{i=-n}^{-1} \tau_i\right) > x\right) = o(\overline{F}^s(x)).$$

(See [4] for the proof in the single-server queue case.)

4.2. *Tail asymptotics for the supremum of a random walk with subexponential increments.* We now remind the well-known result from [14] and [22]. We use negative indices in order to link the result with queueing applications.

THEOREM 3. *Let $\{\xi_n\}$ be an i.i.d. sequence with negative mean $\mathbf{E}\xi_1 = -\alpha$, $S_0 = 0$, $S_{-n} = \sum_1^n \xi_{-i}$ and $\overline{S} = \sup_{n \geq 0} S_{-n}$. Assume that there exists a distribution function $F$ on $[0, \infty)$ such that $F^s$ is subexponential and $\mathbf{P}(\xi_1 > x) \sim d\overline{F}(x)$ with $d > 0$ as $x \to \infty$. Then, as $x \to \infty$,*

$$\mathbf{P}(\overline{S} > x) = (1 + o(1))\frac{d}{\alpha}\overline{F}^s(x).$$

*In particular, consider a $GI/GI/1/\infty$ queue with i.i.d. service times $\{\sigma_n\}$ (with mean $b$) and i.i.d. interarrival times $\{\tau_n\}$ (with mean $a > b$) and put $\xi_n = \sigma_n - \tau_n$. Assume that $\mathbf{P}(\sigma_1 > x) \sim d\overline{F}(x)$, with $F$ as above. Then the stationary waiting time $W$ and the stationary response time $R$ of customer $0$ are such that*

$$\mathbf{P}(R > x) \sim \mathbf{P}(W > x) = (1 + o(1))\frac{d}{a - b}\overline{F}^s(x).$$

*In particular, if the distribution function of $\sigma$ is $F$, then $\mathbf{P}(R > x) \sim \mathbf{P}(W > x) = (1 + o(1))\frac{1}{a-b}\overline{F}^s(x)$.*

The following lower bound is also known (see, e.g., [4]) and was obtained by the use of the strong law of large numbers (SLLN):

THEOREM 4. *Consider a $G/G/1/\infty$ queue with i.i.d. service times $\{\sigma_n\}$ (with mean $b$) and independent stationary ergodic interarrival times $\{\tau_n\}$ (with mean $a > b$). Assume that $\mathbf{P}(\sigma_1 > x) \sim d\overline{F}(x)$, where $d \geq 0$ and where the integrated distribution $F^s$ is long tailed. Then*

$$\mathbf{P}(R > x) \geq \mathbf{P}(W > x) \geq \frac{d}{a - b}\overline{F}^s(x) + o(\overline{F}^s(x)).$$

4.3. *Bounds.*

4.3.1. *Upper bound.* Let $Z = Z_{(-\infty, 0]}$ and let $L$ be the integer defined in Section 2.5 and let $\widehat{s}$ be the service time in the associated $L$-upper-bound $GI/GI/1/\infty$ queue.

Put $\widehat{b} = \mathbf{E}\widehat{s}$ and note that $\mathbf{E}\widehat{\tau} = La$. Then $\widehat{\rho} = \frac{\widehat{b}}{La} = \lambda\gamma(0)(1+\delta) < 1$ where $\delta$ may be chosen as small as possible. We deduce from (7) and (H) that

$$\mathbf{P}(\widehat{s}_1 > x) \sim dL\overline{F}(x).$$

Thus, from Theorem 3,

$$\mathbf{P}(\widehat{R} > x) \sim \frac{1}{La - \widehat{b}} \int_x^\infty \mathbf{P}(\widehat{s} > y) \, dy$$

$$\sim \frac{1}{La - \widehat{b}} \int_x^\infty dL\overline{F}(y) \, dy = \frac{d}{a - \widehat{b}/L} \overline{F}^s(x).$$

Here $\widehat{b}/L \to \gamma(0)$ as $L \to \infty$. We have proved the following.

THEOREM 5. *Under the* (IA), (AA), (SE) *and* (H) *assumptions,*

$$(21) \qquad \limsup_{x \to \infty} \frac{\mathbf{P}(Z > x)}{\overline{F}^s(x)} \leq \lim_{L \to \infty} \lim_{x \to \infty} \frac{\mathbf{P}(\widehat{R} > x)}{\overline{F}^s(x)} = \frac{d}{a - \gamma(0)}.$$

REMARK 4. The assumptions of Theorem 5 bear on the random variables $Y_1^{(j)}$. These can be weakened by considering conditions on the random variables $Z_n = Z_{[n,n]}$ as follows: If the random variables $Z_n$ are i.i.d. with distribution $G$ such that both $G$ and $G^s$ are subexponential, and if the random variables $\tau_n$ are i.i.d. and independent of the $\{Z_n\}$ sequence, then

$$(22) \qquad \limsup_{x \to \infty} \frac{\mathbf{P}(Z > x)}{\overline{G}^s(x)} \leq \frac{1}{a - \gamma(0)}.$$

The proof of this is based on Corollary 2 and on coupling arguments.

4.3.2. *Lower bound.* From (9),

$$Z = Z_{(-\infty, 0]} \geq \underline{R} = \max_j \sup_{n \leq 0} \left( \sum_{i=n}^0 Y_i^{(j)} - \sum_{i=n}^{-1} \tau_i \right) \equiv \max_j R^{(j)}.$$

Then from Theorem 3,

$$\mathbf{P}(R^{(j)} > x) \sim \frac{d^{(j)}}{a - b^{(j)}} \overline{F}^s(x),$$

with $b^{(j)} = \mathbf{E} Y_1^{(j)}$. Note that

$$(23) \qquad \begin{aligned} \sum_j \mathbf{P}(R^{(j)} > x) &\geq \mathbf{P}\left(\max_j R^{(j)} > x\right) \\ &\geq \sum_j \mathbf{P}(R^{(j)} > x) - \sum_{i_1 \neq i_2} \mathbf{P}(R^{(i_1)} > x, R^{(i_2)} > x). \end{aligned}$$

Since, for any $i_1 \neq i_2$,

$$(24) \qquad \mathbf{P}(R^{(i_1)} > x, R^{(i_2)} > x) = o(\overline{F}^s(x))$$

(see Section A.1 for the proof), we get

$$\mathbf{P}(\underline{R} > x) = \sum_j \mathbf{P}(R^{(j)} > x) + o(\overline{F}^s(x)).$$

Thus, the following theorem holds.

THEOREM 6. *Under Assumptions* (IA), (AA), (SE) *and* (H),

$$(25) \qquad \liminf_{x \to \infty} \frac{\mathbf{P}(Z > x)}{\overline{F}^s(x)} \geq \lim_{x \to \infty} \frac{\mathbf{P}(\underline{R} > x)}{\overline{F}^s(x)} = \sum_{j=1}^r \frac{d^{(j)}}{a - b^{(j)}}.$$

REMARK 5. The asymptotics for the lower and upper bounds are the same up to multiplicative constants. So Theorems 5 and 6 imply $\mathbf{P}(Z > x) = \Theta(\overline{F}^s(x))$.

In the single-server isolated queue case, $\gamma(0) = b = M$, $\widehat{b} = Lb$ and $d = 1$. Therefore, in this case the upper and lower bounds coincide.

### 4.4. Examples.

#### 4.4.1. Tandem queues.
The definitions and notation are those of Section 2.6.1. We assume that

$$(26) \qquad \overline{F}_i(x) = \mathbf{P}(\sigma^{(i)} > x) \sim d^{(i)} \overline{F}(x),$$

that $d \equiv d^{(1)} + d^{(2)} > 0$ and that both $F$ and $F^s$ are subexponential. Assumption (H) is valid if we assume in addition that $\sigma_n^{(1)}$ and $\sigma_n^{(2)}$ are independent.

Denote by $Z$ the stationary sojourn time in the network. We look for the asymptotic behavior of the function $\mathbf{P}(Z > x)$ as $x \to \infty$.

The lower bound (25) is

$$\liminf_{x \to \infty} \frac{\mathbf{P}(Z > x)}{\overline{F}^s(x)} \geq \frac{d^{(1)}}{a - b^{(1)}} + \frac{d^{(2)}}{a - b^{(2)}}.$$

Since $\gamma(0) = b \equiv \max(b^{(1)}, b^{(2)})$, the upper bound (21) reads

$$\limsup_{x \to \infty} \frac{\mathbf{P}(Z > x)}{\overline{F}^s(x)} \leq \frac{d}{a - b}.$$

This upper bound was proved in [10].

#### 4.4.2. Generalized Jackson networks.
The definitions are notation and those of Section 2.6.3 and of Section A.2 in the Appendix. We assume that $\mathbf{P}(\sigma^{(i)} > x) \sim l^{(j)} \overline{F}(x)$, that $\sum_j l^{(j)} > 0$ and that both $F$ and $F^s$ are subexponential. Put $\pi^{(j)} = \mathbf{E}\nu^{(j)}$ and $b^{(j)} = \mathbf{E} Y^{(j)} \equiv \pi^{(j)} \mathbf{E}\sigma^{(j)}$. Without loss of generality we may assume $b^{(j)}$ to be positive for all $j$. The network is stable if $\gamma(0) < a$, where $\gamma(0) = b \equiv \max(b^{(1)}, \ldots, b^{(r)})$ (see, e.g., [5]). Assumption (H) is valid (see the example at the end of Section A.2).

Denote by $Z$ the stationary maximal dater. Then, from (25), the lower bound for $\mathbf{P}(Z > x)$ is

$$\liminf_{x \to \infty} \frac{\mathbf{P}(Z > x)}{\overline{F}^s(x)} \geq \sum_1^r \frac{l^{(j)} \pi^{(j)}}{a - b^{(j)}}$$

and the upper bound (21) reads

$$\limsup_{x \to \infty} \frac{\mathbf{P}(Z > x)}{\overline{F}^s(x)} \leq \frac{\sum l^{(j)} \pi^{(j)}}{a - b}.$$

4.4.3. *Max-plus networks.* Similar bounds were studied within the framework of open, irreducible max-plus networks in [10]. As in the single-server isolated queue case (which is an instance of such networks), the upper and lower bounds coincide, which yields the exact asymptotics. However, the exact asymptotics are not known for reducible max-plus networks, even for the particular case of tandem queues.

## 5. Typical event of a subexponential monotone-separable network.

5.1. *Typical event of a subexponential $GI/GI/1$ queue.* This section contains qualitative indications on how rare events occur in subexponential $GI/GI/1$ queues in terms of asymptotic equivalences involving the so-called typical event.

Results of the same nature were first stated by Anantharam in [1] in the regularly varying case (see Theorem 3.1 therein) and by Asmussen and Klüppelberg in [2]. However, we could not find the equivalences (Corollary 5) precisely needed for the extension to monotone separable networks (Theorem 8) in any earlier paper. Some notation and ideas of the proof of Corollary 5 will be used in Theorem 8.

Consider a $GI/GI/1/\infty$ queue with mean interarrival times $a = \mathbf{E}\tau_n$ and mean service times $b = \mathbf{E}\sigma_n$, where $a > b$. Denote by $F$ the distribution function of $\sigma$ and assume that $F$ satisfies (b) of (SE) (i.e., $F^s$ is subexponential), and let $N_x$ be the associated function defined in Section 4.1.2. Let

$$\xi_n = \sigma_n - \tau_n, \qquad S_n^\tau = \sum_1^n \tau_{-i}, \qquad S_n^\sigma = \sum_1^n \sigma_{-i}, \qquad S_n = \sum_1^n \xi_{-i} \equiv S_n^\sigma - S_n^\tau.$$

COROLLARY 5. *Let $W$ (resp. $R$) denote the stationary waiting (resp. response) time of customer 0 in the FIFO $GI/GI/1/\infty$ queue. For any $x$, let $\{K_{n,x}\}$ be a sequence of events such that:*

(i) For any $n$, the event $K_{n,x}$ and the random variable $\sigma_{-n}$ are independent;

(ii) $\mathbf{P}(K_{n,x}) \to 1$ uniformly in $n \geq N_x$ as $x \to \infty$.

For any sequence $\eta_n \to 0$, let

$$(27) \quad A_{n,x} = K_{n,x} \cap \{\sigma_{-n} > x + n(a - b + \eta_n)\} \quad \text{and} \quad A_x = \bigcup_{n \geq N_x} A_{n,x}.$$

Then, as $x \to \infty$,

$$(28) \quad \begin{aligned} \mathbf{P}(W > x) &\sim \mathbf{P}(W > x, A_x) \sim \mathbf{P}(A_x) \\ &\sim \sum_{n=N_x}^{\infty} \mathbf{P}(W > x, A_{n,x}) \sim \sum_{n=N_x}^{\infty} \mathbf{P}(A_{n,x}) \end{aligned}$$

and

$$(29) \quad \mathbf{P}(R > x) \sim \mathbf{P}(W > x).$$

PROOF. Simple calculations using the fact that $F^s$ is long tailed show that, as $x \to \infty$,

$$\sum_{n \geq N_x} \mathbf{P}(A_{n,x}) = \sum_{n \geq N_x} \mathbf{P}(K_{n,x}) \mathbf{P}(\sigma_{-n} > x + n(a - b + \eta_n))$$

$$\sim \sum_{n \geq N_x} \overline{F}(x + n(a - b + \eta_n)) \sim \frac{1}{a - b} \overline{F}^s(x).$$

Thus, if the sequences $\{K_{n,x}\}$ and $\{\eta_n\}$ are such that, for all sufficiently large $x$:

(a) the events $A_{n,x}$ are disjoint for all $n \geq N_x$;
(b) $A_{n,x} \subseteq \{W > x\}$ for all $n \geq N_x$;

then

$$\mathbf{P}(W > x) \geq \mathbf{P}(W > x, A_x) = \mathbf{P}(A_x)$$

$$= \sum_{n \geq N_x} \mathbf{P}(W > x, A_{n,x}) = \sum_{n \geq N_x} \mathbf{P}(A_{n,x}) \sim \frac{1}{a - b} \overline{F}^s(x).$$

Combining this with Theorem 3, we get the equivalences (28).

We now construct two specific sequences $\{K_{n,x}\}$ and $\{\eta_n\}$ satisfying (a), (b) and the assumptions of the corollary. Due to the SLLN, there exists a nonincreasing sequence $\varepsilon_n \to 0$ such that $n\varepsilon_n \to \infty$ and, as $n \to \infty$,

$$\mathbf{P}\left(\left|\frac{S_k^\tau}{k} - a\right| \leq \varepsilon_k, \left|\frac{S_k^\sigma}{k} - b\right| \leq \varepsilon_{k+1} \ \forall k \geq n\right) \to 1.$$

Put $\eta_n = 3\varepsilon_n$ and

$$(30) \quad K_{n,x} = \left\{ \left| \frac{S_k^\tau}{k} - a \right| \leq \varepsilon_k \; \forall \, N_x \leq k \leq n \right\} \cap \left\{ \left| \frac{S_k^\sigma}{k} - b \right| \leq \varepsilon_{k+1} \; \forall \, N_x \leq k < n \right\}.$$

Clearly, the conditions of the corollary are satisfied. Since $n\varepsilon_n > b$ for all sufficiently large $n$, on the event $A_{n,x}$,

$$W \geq S_n > x + n\eta_n - (2n-1)\varepsilon_n - b \geq x.$$

In addition, the events $A_{n,x}$, $n \geq N_x$, are disjoint if $\varepsilon_{N_x} \leq (a-b)/2$. Indeed, on the event $A_{n,x}$, we then have $S_n > x$ and $S_{n-1}^* = \max_{0 \leq j \leq n-1} S_j \leq \max_{0 \leq j \leq n-1} j(b - a + 2\varepsilon_{N_x}) \leq 0$; and the events $\{S_{n-1}^* \leq 0\} \cap \{S_n > x\}$ are obviously disjoint.

Take now any other sequences $\{\tilde{K}_{n,x}\}$ and $\{\tilde{\eta}_n\}$ satisfying the conditions of the corollary and denote the corresponding events by $\{\tilde{A}_{n,x}\}$ and $\tilde{A}_x$. Then

$$|\mathbf{P}(A_x) - \mathbf{P}(\tilde{A}_x)|$$
$$\leq \mathbf{P}\left( \bigcup_{n \geq N_x} \{\sigma_{-n} > x + n(a - b + \min(\eta_n, \tilde{\eta}_n))\} \right)$$
$$\quad - \mathbf{P}\left( \bigcup_{n \geq N_x} K_{n,x} \cap \tilde{K}_{n,x} \cap \{\sigma_{-n} > x + n(a - b + \max(\eta_n, \tilde{\eta}_n))\} \right)$$
$$\leq \sum_{n \geq N_x} \mathbf{P}(\sigma_{-n} > x + n(a - b + \min(\eta_n, \tilde{\eta}_n)))$$
$$\quad - \sum_{n \geq N_x} \mathbf{P}(K_{n,x} \cap \tilde{K}_{n,x} \cap \{\sigma_{-n} > x + n(a - b + \max(\eta_n, \tilde{\eta}_n))\})$$
$$\leq \Delta_x \sum_{n \geq N_x} \overline{F}(x + n(a - b + \max(\eta_n, \tilde{\eta}_n)))$$
$$\quad + \sum_{n \geq N_x} \Big( \overline{F}(x + n(a - b + \min(\eta_n, \tilde{\eta}_n)))$$
$$\quad\quad\quad - \overline{F}(x + n(a - b + \max(\eta_n, \tilde{\eta}_n))) \Big),$$

where $\Delta_x = \sup_{n \geq N_x}(\mathbf{P}(K_{n,x}^c) + \mathbf{P}(\widetilde{K}_{n,x}^c)) \to 0$ as $x \to \infty$. Thus, both terms in the last expression are $o(\overline{F}^s(x))$, and the equivalences (28) hold for sequences $\{\widetilde{K}_{n,x}\}$ and $\{\tilde{\eta}_n\}$.

Finally, the equivalence (29) follows from the relation $R = W + \sigma_0$, the independence and the fact that the tail of $W$ is heavier than that of $\sigma$ [see (15)]. □

The event $A_x$ (which will be referred to as the *typical event* of the subexponential $GI/GI/1$ queue in what follows) occurs if there is only one big service time and all other service times or interarrival times follow the SLLN.

5.2. *Key equivalences for the maximal dater of a subexponential monotone-separable network.* In this section, we consider a monotone-separable network satisfying (IA), (AA), (H) and (SE). The function $N_x$ is that associated with the reference distribution function $F$ of the (SE) assumptions.

THEOREM 7. *Let $Z$ be the stationary maximal dater of some monotone separable network. Denoting $\widehat{A}_x$ the typical event of the L-upper-bound queue (more generally we will add a hat to indicate that a variable pertains to the upper bound queue), we have*

$$(31) \qquad \mathbf{P}(Z > x) \sim \mathbf{P}(Z > x, \widehat{A}_x) \sim \sum_{n=N_x}^{\infty} \mathbf{P}(Z > x, \widehat{A}_{n,x})$$

*and*

$$(32) \qquad \mathbf{P}(Z > x, \widehat{A}_x) = \Theta(\overline{F}^s(x)).$$

*Also, for any random variable $\tilde{Z}$ such that $\tilde{Z} \leq Z$ a.s.,*

$$(33) \qquad \mathbf{P}(\tilde{Z} > x) = \sum_{n \geq N_x} \mathbf{P}(\tilde{Z} > x, \widehat{A}_{n,x}) + o(\overline{F}^s(x)).$$

PROOF. Since $Z \leq \widehat{R}$ a.s.,

$$\begin{aligned}\mathbf{P}(Z > x) &= \mathbf{P}(Z > x, \widehat{A}_x) + \mathbf{P}(Z > x; \widehat{R} > x, (\widehat{A}_x)^c) \\ &= \mathbf{P}(Z > x, \widehat{A}_x) + o(\overline{F}^s(x)) = \sum_{n \geq N_x} \mathbf{P}(Z > x, \widehat{A}_{n,x}) + o(\overline{F}^s(x))\end{aligned}$$

from Corollary 5.

From Theorems 5 and 6, $\mathbf{P}(Z > x) = \Theta(\overline{F}^s(x))$. Thus,

$$\mathbf{P}(Z > x, \widehat{A}_x) = \mathbf{P}(Z > x) - o(\overline{F}^s(x)) = \Theta(\overline{F}^s(x)) - o(\overline{F}^s(x)) = \Theta(\overline{F}^s(x))$$

and both (31) and (32) follow. □

The main result of the paper concerning subexponential monotone separable networks is the following theorem, which can be seen as a network extension of Corollary 5 and which gives the shape of the typical event creating a large maximal dater in such a network.

THEOREM 8. *The assumptions are the same as in Theorem 7. Put $b = \gamma(0)$. For any $x$ and for $j = 1, \ldots, r$, let $\{K_{n,x}^{(j)}\}$ be a sequence of events such that:*

(i) *For any $n$ and $j$, the event $K_{n,x}^{(j)}$ and the random variable $Y_{-n}^{(j)}$ are independent;*

(ii) *For any $j$, $\mathbf{P}(K_{n,x}^{(j)}) \to 1$ uniformly in $n \geq N_x$ as $x \to \infty$.*

*For all sequences $\eta_n^{(j)}$, $j = 1, \ldots, r$, tending to 0, put*

(34)
$$A_{n,x}^{(j)} = K_{n,x}^{(j)} \cap \{Y_{-n}^{(j)} > x + n(a - b + \eta_n^{(j)})\},$$
$$A_x^{(j)} = \bigcup_{n=N_x}^{\infty} A_{n,x}^{(j)} \quad \text{and} \quad A_x = \bigcup_{j=1}^{r} A_x^{(j)}.$$

*Then, as $x \to \infty$,*

(35)
$$\mathbf{P}(Z > x) \sim \mathbf{P}(Z > x, A_x) \sim \sum_{1}^{r} \mathbf{P}(Z > x, A_x^{(j)})$$
$$\sim \sum_{j=1}^{r} \sum_{n=N_x}^{\infty} \mathbf{P}(Z > x, A_{n,x}^{(j)}).$$

*Similarly, for any random variable s.t. $\tilde{Z} \leq Z$,*

(36)
$$\mathbf{P}(\tilde{Z} > x) = \mathbf{P}(\tilde{Z} > x, A_x) + o(\overline{F}^s(x))$$
$$= \sum_{1}^{r} \mathbf{P}(\tilde{Z} > x, A_x^{(j)}) + o(\overline{F}^s(x))$$
$$= \sum_{j=1}^{r} \sum_{n=N_x}^{\infty} \mathbf{P}(\tilde{Z} > x, A_{n,x}^{(j)}) + o(\overline{F}^s(x)).$$

*If $\mathbf{P}(\tilde{Z} > x) = \Theta(\overline{F}^s(x))$, one can replace the last equalities by equivalences and delete the $o(\overline{F}^s(x))$ terms in the last relation.*

The equivalences (31) and (35) will be the key relationships for the exact asymptotics of the examples of Section 6. They show that for the monotone separable network also, whenever the maximal dater is large, at most one of the service times is large whereas all other ones are moderate.

PROOF. We will only prove the equivalence

(37)
$$\mathbf{P}(Z > x) \sim \sum_{j=1}^{r} \sum_{n=N_x}^{\infty} \mathbf{P}(Z > x, A_{n,x}^{(j)}).$$

The other equivalences in (35) may be obtained similarly.

We start the proof with the following three reductions.

First, it is sufficient to prove the equivalence (37) when replacing each of the events $K_{n,x}^{(j)}$ by the whole probability space $\Omega$. Indeed, put $\tilde{A}_{n,x}^{(j)} = \{Y_{-n}^{(j)} >$ $x + n(a - b + \eta_n^{(j)})\}$. We know from Theorem 7 that $\mathbf{P}(Z > x) = \Theta(\overline{F}^s(x))$. Suppose that

$$\mathbf{P}(Z > x) \sim \sum_{j=1}^{r} \sum_{n \geq N_x} \mathbf{P}(Z > x, \tilde{A}_{n,x}^{(j)}).$$

Then

$$\sum_{j=1}^{r} \sum_{n \geq N_x} \mathbf{P}(Z > x, A_{n,x}^{(j)})$$

$$= \sum_{j=1}^{r} \sum_{n \geq N_x} \mathbf{P}(Z > x, \tilde{A}_{n,x}^{(j)}) - \sum_{j=1}^{r} \sum_{n \geq N_x} \mathbf{P}(Z > x, \tilde{A}_{n,x}^{(j)} \setminus A_{n,x}^{(j)}).$$

The result then follows from the fact that the last subtracted sum is non-negative and is not bigger than

$$\sum \sum \mathbf{P}(\tilde{A}_{n,x}^{(j)}) \mathbf{P}((K_{n,x}^{(j)})^c) \leq \Delta(x) \Theta(\overline{F}^s(x)) = o(\overline{F}^s(x))$$

since $\Delta(x) \equiv \max_{1 \leq j \leq r} \sup_{n \geq N_x} \mathbf{P}((K_{n,x}^{(j)})^c) \to 0$ as $x \to \infty$.

Second, it is sufficient to consider the case $\eta_n^{(j)} = 0$ for all $n$ and $j$. This follows from the following bound where $\delta_x = \max_j \sup_{n \geq N_x} (\eta_n^{(j)})^+$:

$$\sum_{j=1}^{r} \sum_{n \geq N_x} \mathbf{P}(Z > x, Y_{-n}^{(j)} > x + n(a - b))$$
$$- \mathbf{P}(Z > x, Y_{-n}^{(j)} > x + n(a - b + (\eta_n^{(j)})^+))$$
$$\leq \sum_{j=1}^{r} \sum_{n \geq N_x} \mathbf{P}(Y_{-n}^{(j)} \in (x + n(a - b), x + n(a - b + \delta_x)))$$
$$= (1 + o(1))\, d\overline{F}^s(x) \left( \frac{1}{a-b} - \frac{1}{a-b+\delta_x} \right) = o(\overline{F}^s(x))$$

and a symmetrical bound for the negative part of $\eta_n^{(j)}$. Thus it is enough to prove the equivalence

(38) $$\mathbf{P}(Z > x) \sim \sum_{j=1}^{r} \sum_{n \geq N_x} \mathbf{P}(Z > x, Y_{-n}^{(j)} > x + n(a - b)).$$

Third, if there exists a sequence $\varepsilon_L \in (0, a - b)$, $\varepsilon_L \to 0$ such that, for any $L$, the following equivalence takes place (where $b_L = b + \varepsilon_L$):

(39) $$\mathbf{P}(Z > x) \sim \sum_{j=1}^{r} \sum_{n \geq N_x} \mathbf{P}(Z > x, Y_{-n}^{(j)} > x + n(a - b_L)),$$

then (38) holds. Indeed, take $\varepsilon_L < (a-b)/2$. Then

$$\sum_{j=1}^{r} \sum_{n \geq N_x} \mathbf{P}(Z > x, Y^{(j)}_{-n} > x + n(a - b_L))$$

$$- \mathbf{P}(Z > x, Y^{(j)}_{-n} > x + n(a - b))$$

$$\leq \sum_{j=1}^{r} \sum_{n \geq N_x} \mathbf{P}(Y^{(j)}_{-n} \in [x + n(a - b_L), x + n(a - b)])$$

$$= (1 + o(1)) d\overline{F}^s(x) \left( \frac{1}{a - b_L} - \frac{1}{a - b} \right)$$

$$\leq (1 + o(1)) \frac{2\varepsilon_L d}{(a - b)^2} \overline{F}^s(x).$$

Letting $L \to \infty$, we derive (38) from (39).

Before proving (39), we recall that, from conditions (SE) and (H),

$$\mathbf{P}(\widehat{s}_1 > x) \sim \mathbf{P}\left( \bigcup_{j=1}^{r} \bigcup_{l=1}^{L} \{Y^{(j)}_l > x\} \right) \sim \sum_{j=1}^{r} \sum_{l=1}^{L} \mathbf{P}(Y^{(j)}_l > x).$$

Since $F$ is long tailed, we can replace the latter equivalences by

(40)
$$\mathbf{P}(\widehat{s}_1 > x) \sim \mathbf{P}\left( \bigcup_{j=1}^{r} \bigcup_{l=1}^{L} \{Y^{(j)}_l > x + l(a - b_L)\} \right)$$
$$\sim \sum_{j=1}^{r} \sum_{l=1}^{L} \mathbf{P}(Y^{(j)}_l > x + l(a - b_L)).$$

More precisely, when denoting the event in the left-hand side by $C_x$ and the event in the center by $D_x$, we get $C_x \subseteq D_x$ and

(41)
$$\sup_{y \geq x} \frac{\mathbf{P}(D_y \setminus C_y)}{\overline{F}(y)} = o(1),$$

when $x \to \infty$, while $\mathbf{P}(C_x) = \Theta(\overline{F}(x))$.

We now prove (39). For any $L$, put $\tilde{N}_x = N_x/L$ (more precisely the integer part of this ratio) and note that $\tilde{N}_x$ also satisfies condition (16). Take $L$ sufficiently large and set $b_L = \mathbf{E}\,\widehat{s}_0/L \equiv b + \varepsilon_L$. For the $L$-upper queue, one can take the typical event of the form

$$\widehat{A}_x = \bigcup_{n \geq \tilde{N}_x} \{\widehat{s}_{-n} > x + nL(a - b_L)\} \equiv \bigcup_{n \geq \tilde{N}_x} \widehat{A}_{n,x}.$$

From Theorem 7,

$$\mathbf{P}(Z > x) \sim \sum_{m \geq \tilde{N}_x} \mathbf{P}(Z > x, \widehat{s}_{-m} > x + mL(a - b_L)).$$

From (40) and (41),

$$\sum_{m \geq \tilde{N}_x} \mathbf{P}(Z > x, \hat{s}_{-m} > x + mL(a - b_L))$$

$$= (1 + o(1))$$

$$\times \sum_{m \geq \tilde{N}_x} \sum_{j=1}^{r} \sum_{l=1}^{L} \mathbf{P}(Z > x, Y^{(j)}_{-mL+l} > x + mL(a - b_L) + l(a - b_L)),$$

where the uniformity in $m$ required to obtain the term $o(1)$ follows from the uniformity in $y$ in (41). So

$$\mathbf{P}(Z > x) \sim \sum_{j=1}^{r} \sum_{n \geq N_x} \mathbf{P}(Z > x, Y^{(j)}_{-n+1} > x + n(a - b_L))$$

$$\sim \sum_{j=1}^{r} \sum_{n \geq N_x} \mathbf{P}(Z > x, Y^{(j)}_{-n} > x + n(a - b_L)). \qquad \square$$

**6. Two examples of exact tail asymptotics.** This section gives two illustrations of the use of Theorem 8 in order to derive exact asymptotics. Without loss of generality we can assume interarrival times to be constants and equal to $a$ (see Section A.3).

6.1. *Tandem queues.* For tandem queue, the assumptions on the tails of the service times are those of Section 4.4.1. Here $Y^{(j)}_n = \sigma^{(j)}_n$, so that (H) trivially holds, since the service times are independent. The results are stated for the two-station case, but the extension to tandems (or treelike networks) of any dimension is immediate.

Choose a sequence $N_x$ satisfying (16). Let $W^{(j)}_n$ be the stationary waiting time of customer $n$ in queue $j = 1, 2$, and $\tau^{(2)}_n$ be the interarrival time between the $n$th and $(n+1)$st customers to the second queue. Let $\xi^{(1)}_n = \sigma^{(1)}_n - \tau_n \equiv \sigma^{(1)}_n - a$, $\xi^{(2)}_n = \sigma^{(2)}_n - \tau^{(2)}_n$ and $\tilde{\xi}_n = \sigma^{(2)}_n - \sigma^{(1)}_{n+1}$. We will also use the following notation:

$$\tilde{S}_n = \sum_{i=1}^{n} \tilde{\xi}_{-i}, \qquad S^{(j)}_n = \sum_{i=1}^{n} \xi^{(j)}_{-i}, \qquad S^{(\sigma,j)}_n = \sum_{i=0}^{n} \sigma^{(j)}_{-i}, \qquad j = 1, 2.$$

The following relations hold:

(42) $$W^{(j)}_{n+1} = \max(0, W^{(j)}_n + \xi^{(j)}_n), \qquad j = 1, 2,$$

and

(43) $$\tau^{(2)}_n = -\min(0, W^{(1)}_n + \xi^{(1)}_n) + \sigma^{(1)}_{n+1} \geq \sigma^{(1)}_{n+1}$$

so that $\tau_n^{(2)} = \sigma_{n+1}^{(1)}$ if $W_n^{(1)} + \xi_n^{(1)} \geq 0$. In addition,

$$Z = W_0^{(1)} + \sigma_0^{(1)} + W_0^{(2)} + \sigma_0^{(2)}. \tag{44}$$

Also, from (11),

$$Z = \sup_{-\infty < p \leq q \leq 0} \left( \sum_{m=p}^{q} \sigma_m^{(1)} + S_q^{(\sigma,2)} + pa \right). \tag{45}$$

Similarly, for any $n$, the stationary response time $Z_{(-\infty,-n]}$ of customer $(-n)$ satisfies the relations

$$Z_{(-\infty,-n]} = W_{-n}^{(1)} + \sigma_{-n}^{(1)} + W_{-n}^{(2)} + \sigma_{-n}^{(2)}$$

$$= \sup_{-\infty < p \leq q \leq -n} \left( \sum_{m=p}^{q} \sigma_m^{(1)} + \sum_{m=q}^{-n} \sigma_m^{(2)} + (p+n)a \right). \tag{46}$$

6.1.1. *End-to-end delay.* In this section, we prove the following exact asymptotic, which refines the bounds of Section 4.4.1 (these bounds do not coincide in general).

THEOREM 9. *Under the assumptions of Section* 4.4,

$$\mathbf{P}(Z > x) \sim \left( \frac{d^{(1)}}{a - b} + \frac{d^{(2)}}{a - b^{(2)}} \right) \overline{F}^s(x), \tag{47}$$

*where* $b = \max(b^{(1)}, b^{(2)}) \equiv \gamma(0)$.

REMARK 6. As a corollary of Theorem 9 and of results from [3] and [16], one can easily derive sharp asymptotics for the stationary queue length $Q = Q_1 + Q_2$ in the tandem queue. Also, the result may be easily extended to queues in tandem of any finite length.

PROOF OF THEOREM 9. From Theorem 8, we get

$$\mathbf{P}(Z > x) \sim \sum_{j=1}^{2} \mathbf{P}(Z > x, A_x^{(j)}) \sim \sum_{n=N_x}^{\infty} \sum_{j=1}^{2} \mathbf{P}(Z > x, A_{n,x}^{(j)}).$$

We have to find appropriate sequences $\{K_{n,x}^{(j)}\}$ and $\{\eta_n^{(j)}\}$.

Start with $j = 1$. For any $\{K_{n,x}^{(1)}\}$ and $\eta_n^{(1)} \to 0$,

$$\sum_{n \geq N_x} \mathbf{P}(Z > x, A_{n,x}^{(1)}) \leq \sum_{n \geq N_x} \mathbf{P}(\sigma_{-n}^{(1)} > x + n(a-b) + \eta_n^{(1)}) \sim \frac{d^{(1)}}{a - b} \overline{F}^s(x).$$

For the lower bound, consider the events
$$K_{n,x}^{(1)} = \{S_{n-1}^{(\sigma,1)} \geq n(b^{(1)} - \eta_n^{(1)}), S_n^{(\sigma,2)} \geq n(b^{(2)} - \eta_n^{(1)})\}$$
and choose a sequence $\eta_n^{(1)} \to 0$ such that $\mathbf{P}(K_{n,x}^{(1)}) \to 1$ uniformly in $n \geq N_x$ as $x \to \infty$. Then, from (45),

$$\mathbf{P}(Z > x, A_{n,x}^{(1)}) \geq \mathbf{P}\Big(\sigma_{-n}^{(1)} + \max{(S_{n-1}^{(\sigma,1)}, S_n^{(\sigma,2)})} - na > x, A_{x,n}^{(1)}\Big)$$
$$\geq \mathbf{P}\Big(\sigma_{-n}^{(1)} + n(\max{(b^{(1)}, b^{(2)})} - \eta_n^{(1)} - a) > x\Big)\mathbf{P}(K_{n,x}^{(1)})$$
$$= (1 + o(1))\mathbf{P}(\sigma_{-n}^{(1)} > x + n(a - b + \eta_n^{(1)})),$$

and the lower bound for $\mathbf{P}(Z > x, A_x^{(1)})$ is asymptotically equivalent to the upper one.

Consider $j = 2$. The lower bound
$$\mathbf{P}(Z > x, A_x^{(2)}) \geq \mathbf{P}(W_0^{(2)} > x, A_x^{(2)}) = (1 + o(1))\frac{d^{(2)}}{a - b^{(2)}}\overline{F}^s(x)$$
follows from Theorem 4.

For the upper bound, put
$$U_n = \sup_{-\infty < p \leq 0} \sup_{\max(p,-n) < q \leq 0} \left(\sum_{m=p}^{q} \sigma_m^{(1)} + S_q^{(\sigma,2)} + pa\right)$$
and note that, from (45) and (46),
$$Z \leq \max{(Z_{(-\infty,-n-1]} + S_n^{(\sigma,2)} - na, U_n)}$$
$$\equiv \max{(Z_{(-\infty,-n-1]} + S_{n-1}^{(\sigma,2)} + \sigma_{-n}^{(2)} - na, U_n)},$$
where the random vector $(Z_{(-\infty,-n-1]}, U_n, S_{n-1}^{(\sigma,2)})$ is independent of $\sigma_{-n}^{(2)}$.

Since $U_n \leq Z$ a.s., $\mathbf{P}(U_n \leq x) \to 1$ uniformly in $n$ as $x \to \infty$. Since the distribution of $Z_{(-\infty,-n-1]}$ does not depend on $n$, $Z_{(-\infty,-n-1]}/n \to 0$ in probability. Due to the SLLN, $S_{n-1}^{(\sigma,2)}/n \to b^{(2)}$ a.s. Therefore, there exists a sequence $\varepsilon_n \downarrow 0$, $n\varepsilon_n \to \infty$ such that
$$\mathbf{P}(U_n \leq x, Z_{(-\infty,-n-1]} \leq n\varepsilon_n, S_{n-1}^{(\sigma,2)} \leq n(b^{(2)} + \varepsilon_n)) \to 1$$
uniformly in $n \geq N_x$ as $x \to \infty$. Denote the latter event by $K_{n,x}^{(2)}$ and recall that it is independent of $\sigma_{-n}^{(2)}$. Put $\eta_n^{(2)} = -2\varepsilon_n$. Then
$$\mathbf{P}(Z > x, A_{n,x}^{(2)}) \leq \mathbf{P}(\sigma_{-n}^{(2)} + n(b^{(2)} - a) + 2n\varepsilon_n > x, A_{n,x}^{(2)})$$
$$= \mathbf{P}(\sigma_{-n}^{(2)} > x + n(a - b^{(2)} + \eta_n^{(2)}), K_{n,x}^{(2)})$$
$$= (1 + o(1))\mathbf{P}(\sigma_{-n}^{(2)} > x + n(a - b^{(2)} + \eta_n^{(2)})),$$
and the desired asymptotics follow. □

6.1.2. *Delay at the second queue.* In this section, we focus on the asymptotics for the stationary waiting time $W^{(2)} \equiv W_0^{(2)}$ of customer 0 at the second queue. The assumptions are the same as in Section 6.1.

Results on the matter were obtained by Huang and Sigman in [19] in the case where the tail of $\sigma^{(2)}$ is heavier than that of $\sigma^{(1)}$. The results of the present section are more general in that such an assumption is not required.

First, let us see how the results of [19] follow from what we have here. Under the assumptions of Section 4.4.1, we get from (36) of Theorem 8 (for $\widetilde{Z} = W^{(2)} \leq Z$) that

(48) $\mathbf{P}(W^{(2)} > x) = \mathbf{P}(W^{(2)} > x, A_x^{(1)}) + \mathbf{P}(W^{(2)} > x, A_x^{(2)}) + o(\overline{F}^s(x)).$

Then

(49) $$\mathbf{P}(W^{(2)} > x, A_x^{(2)}) = \frac{d^{(2)}}{a - b^{(2)}} \overline{F}^s(x) + o(\overline{F}^s(x)).$$

The lower bound follows from Theorem 4 and the upper one from the inequality $\mathbf{P}(W^{(2)} > x, A_x^{(2)}) \leq \mathbf{P}(Z > x, A_x^{(2)})$ and from the part $j = 2$ of the proof of Theorem 9. For the first term in the right-hand side of (48), we have

$$0 \leq \mathbf{P}(W^{(2)} > x, A_x^{(1)}) \leq \mathbf{P}(Z > x, A_x^{(1)}) = \frac{d^{(1)}}{a - b} \overline{F}^s(x) + o(\overline{F}^s(x))$$

from the proof of Theorem 9. Thus, if $d^{(1)} = 0 < d^{(2)}$, then

(50) $$\mathbf{P}(W^{(2)} > x) \sim \frac{d^{(2)}}{a - b^{(2)}} \overline{F}^s(x)$$

which is the result of [19].

We now successively consider the three cases $b^{(1)} > b^{(2)}$, $b^{(1)} = b^{(2)}$ and $b^{(1)} < b^{(2)}$.

*Case $b^{(1)} > b^{(2)}$.* For the following theorem, we do not need any assumption on the tail of $\sigma^{(1)}$. In fact, we do not even need to assume that $F$ is subexponential, the only required assumption being that $F^s$ is subexponential.

THEOREM 10. *Assume $a > b^{(1)} > b^{(2)}$ and $\mathbf{P}(\sigma^{(2)} > x) \sim d^{(2)} \overline{F}(x)$ where $d^{(2)} > 0$, and the integrated tail distribution $F^s$ is subexponential. Then, as $x \to \infty$,*

(51) $$\mathbf{P}(W^{(2)} > x) \sim \frac{d^{(2)}}{a - b^{(2)}} \overline{F}^s(x).$$

PROOF. We already established the right lower bound in (49). Thus, it is enough to derive an upper bound which is asymptotically equivalent to the lower one.

For this, we use the notation from the beginning of Section 6.1. Since $\sigma_i^{(2)}$ and $\tau_i^{(2)}$ are independent and $F^s$ is long tailed, as $x \to \infty$,

$$\text{(52)} \qquad \int_x^\infty \mathbf{P}(\xi_i^{(2)} > t)\, dt \sim \int_x^\infty \mathbf{P}(\sigma_i^{(2)} > t)\, dt \sim d^{(2)} \overline{F}^s(x).$$

Put $S^{(2)} = \sup_{n \geq 0} S_n^{(2)}$ and $\widetilde{S} = \sup_{n \geq 0} \widetilde{S}_n$. Since $\tau_i^{(2)} \geq \sigma_{i+1}^{(1)}$, $W^{(2)} = S^{(2)} \leq \widetilde{S}$ a.s. Since $b^{(1)} > b^{(2)}$, $\widetilde{S} < \infty$ a.s. So, we have

$$\text{(53)} \qquad \begin{aligned} \mathbf{P}(S^{(2)} > x) &= \mathbf{P}(S^{(2)} > x, \widetilde{S} > x) \\ &\leq \sum_n \mathbf{P}(S^{(2)} > x, \widetilde{\xi}_{-n} > x + n\widetilde{c}) + o(F^s(x)), \end{aligned}$$

where (53) follows from Corollary 5 which implies that

$$\{\widetilde{S} > x\} = \bigcup_n \{\widetilde{\xi}_{-n} > x + n\widetilde{c}\} \cup B_x \qquad \text{where } \mathbf{P}(B_x) = o(\overline{F}^s(x)).$$

Set $c = a - b^{(2)}$ and $\widetilde{c} = b^{(1)} - b^{(2)}$. For all $\varepsilon \in (0, \widetilde{c})$, $R > 0$ and $n$, define the event

$$D_{n,\varepsilon,R} = \{S_i^{(2)} \leq R - i(c - \varepsilon), \widetilde{S}_i \leq R - i(\widetilde{c} - \varepsilon),\ i = 1, 2, \ldots, n-1;$$
$$\widetilde{S}_{n+j} - \widetilde{S}_n \leq R - j(\widetilde{c} - \varepsilon),\ j = 1, 2, \ldots\}.$$

By the SLLN, for any $\varepsilon > 0$, there exists $R > 0$ such that, for any $n = 1, 2, \ldots$, $\mathbf{P}(D_{n,\varepsilon,R}) \geq 1 - \varepsilon$. From (53), we have

$$\text{(54)} \qquad \begin{aligned} \mathbf{P}(S^{(2)} > x) &\leq \sum_n \mathbf{P}(S^{(2)} > x, \widetilde{\xi}_{-n} > x + n\widetilde{c}) + o(F^s(x)) \\ &\leq \sum_n \mathbf{P}((D_{n,\varepsilon,R})^c, \widetilde{\xi}_{-n} > x + n\widetilde{c}) \\ &\quad + \sum_n \mathbf{P}(D_{n,\varepsilon,R}, S^{(2)} > x) + o(F^s(x)) \\ &\equiv \Sigma_1 + \Sigma_2 + o(F^s(x)). \end{aligned}$$

Also,

$$\Sigma_1 \leq (1 + o(1))\varepsilon \sum_n \mathbf{P}(\widetilde{\xi}_{-n} > x + n\widetilde{c}) = (1 + o(1))\frac{\varepsilon d^{(2)}}{\widetilde{c}} \overline{F}^s(x).$$

On the intersection of the events $D_{n,\varepsilon,R}$ and $\{\xi_n^{(2)} \leq x - 2R + (n-1)(c - \varepsilon)\}$, we have $S_n^{(2)} \leq x - R$. In addition, $S_i^{(2)} < R$ for $i = 1, \ldots, n-1$ and, for all

$j \geq 1$,
$$S_{n+j}^{(2)} = S_{n+j}^{(2)} + S_n^{(2)} - S_n^{(2)} \leq \widetilde{S}_{n+j} - \widetilde{S}_n + S_n^{(2)} \leq R - j(\widetilde{c} - \varepsilon) + x - R \leq x.$$

Thus, on this intersection, $S_m^{(2)} \leq x$ for all $m$ if $x \geq R$. Therefore,
$$\mathbf{P}(D_{n,\varepsilon,R}, S^{(2)} > x) = \mathbf{P}(D_{n,\varepsilon,R}, S^{(2)} > x, \xi_n^{(2)} > x - 2R + (n-1)(c - \varepsilon))$$
$$\leq \mathbf{P}(\xi_n^{(2)} > x - 2R + (n-1)(c - \varepsilon)).$$

Hence,
$$\Sigma_2 \leq \sum_n \mathbf{P}(\xi_n^{(2)} > x - 2R + (n-1)(c - \varepsilon))$$
$$= (1 + o(1))\frac{d^{(2)}}{c - \varepsilon}\overline{F}^s(x - 2R) = (1 + o(1))\frac{d^{(2)}}{c - \varepsilon}\overline{F}^s(x),$$

as $x \to \infty$, because of (52). Since $\varepsilon > 0$ is arbitrary, the result follows. □

*Case $b^{(1)} = b^{(2)}$.* We assume $v_i^2 = \mathrm{Var}(\sigma^{(i)})$ to be finite for $i = 1, 2$ and we use the notation $v = \sqrt{v_1^2 + v_2^2}$.

THEOREM 11. *Assume $a > b^{(1)} = b^{(2)} \equiv b$ and $\mathbf{P}(\sigma^{(i)} > x) \sim d^{(i)}\overline{F}(x)$ as $x \to \infty$ with $d^{(1)} + d^{(2)} > 0$, where both $F$ and $F^s$ are subexponential. Then, as $x \to \infty$,*

(55)
$$\mathbf{P}(W^{(2)} > x) = 2d^{(1)}\int_0^\infty \overline{F}(x + y(a-b))\overline{\Phi}\left(\frac{x}{v\sqrt{y}}\right)dy$$
$$+ \frac{d^{(2)}}{a - b}\overline{F}^s(x) + o(\overline{F}^s(x)),$$

*where $\overline{\Phi}$ is the tail of the standard normal distribution. In particular, if either:*

(i) $d^{(2)} > 0$ *or*
(ii) $d^{(2)} = 0$, $d^{(1)} > 0$ *and*

(56)
$$\liminf_{x \to \infty} \overline{F}^s(x^2)/\overline{F}^s(x) > 0,$$

*then one can replace the equality in (55) by an equivalence and delete the term $o(\overline{F}^s(x))$ in this equation.*

REMARK 7. Under condition (56), the integral in the right-hand side of (55) is of order $\Theta(\overline{F}^s(x^2)) = \Theta(\overline{F}^s(x))$. Condition (56) is satisfied if the tail $\overline{F}^s(x)$ is "extremely heavy," for example, $\overline{F}^s(x) \sim (\log x)^{-K}$, $K > 0$. However, if $\limsup_{x \to \infty} \overline{F}^s(x^2)/\overline{F}^s(x) = 0$ (the latter holds for Pareto, log-normal and Weibull distributions), then the integral is of order $o(\overline{F}^s(x))$.

PROOF OF THEOREM 11. We use again (48) and (49), and we are left with the problem of finding the asymptotics for

$$\mathbf{P}(W^{(2)} > x, A_x^{(1)}) = \sum_{n \geq N_x} \mathbf{P}(W^{(2)} > x, A_{n,x}^{(1)}) + o(\overline{F}^s(x))$$

for an appropriate event $A_x^{(1)}$ satisfying the assumptions of Theorem 8. Let $S_{m,n} = \sum_{i=1}^m (b - \sigma_{-n+i}^{(1)})$. Due to the LLN, $\max_{1 \leq m \leq n} S_{m,n}/n \to 0$ in probability as $n \to \infty$. Therefore, there exists a nonincreasing sequence $\eta_n^{(1)} \to 0$ such that $\mathbf{P}(\max_{1 \leq m \leq n} S_{m,n}/n > \eta_n^{(1)}) \leq \eta_n^{(1)}$ for all $n$. Take

$$K_{n,x}^{(1)} = \left\{ \max_{1 \leq m \leq n} S_{m,n}/n \leq \eta_n^{(1)} \right\} \cap \{W_{-n-1}^{(2)} + \sigma_{-n-1}^{(2)} \leq x + n(a - b + \eta_n^{(1)})\}$$

$$\equiv K_{n,x}^{(1,1)} \cap K_{n,x}^{(1,2)}.$$

Easy calculations based on the fact that $W_{-n+m}^{(1)} \geq \sigma_{-n}^{(1)} + S_{m,n} - (m+1)a$ show that on the event $A_{n,x}^{(1)} = K_{n,x}^{(1)} \cap \{\sigma_{-n}^{(1)} > x + n(a - b + \eta_n^{(1)})\}$, one has $W_{-n+m}^{(1)} > 0$, for all $m = 1, 2, \ldots, n$. Using the fact that $\tau_{-n-1}^{(2)} \geq \sigma_{-n}^{(1)}$, one gets immediately that on this event, $W_{-n}^{(2)} = 0$. Therefore, on this event, $\tau_j^{(2)} = \sigma_{j+1}^{(1)}$ for all $j = -n+1, \ldots, 0$ and

$$W^{(2)} = W_0^{(2)} = \max_{0 \leq j \leq n} \widetilde{S}_j \equiv V_n.$$

From the central limit theorem for the reflected random walk,

$$\frac{V_n}{v\sqrt{n}} \to \psi$$

weakly, where $\psi$ has the following tail distribution:

$$\mathbf{P}(\psi > x) = 2\overline{\Phi}(x) \equiv \frac{2}{\sqrt{2\pi}} \int_x^\infty \exp\{-y^2/2\} \, dy.$$

Take any $c > 0$. If $N_x \leq n \leq cx^2$, then

$$\mathbf{P}\left(\frac{V_n}{v\sqrt{n}} > \frac{x}{v\sqrt{n}}\right) \leq \mathbf{P}\left(\frac{V_n}{v\sqrt{n}} > \frac{1}{v\sqrt{c}}\right) = (1 + o(1))2\overline{\Phi}\left(\frac{1}{v\sqrt{c}}\right)$$

as $x \to \infty$. For any $\Delta > 0$, choose $c \ll 1$ such that $\overline{\Phi}(\frac{1}{v\sqrt{c}}) \leq \Delta$. Then

$$\sum_{N_x}^{cx^2} \mathbf{P}(W_0^{(2)} > x, A_{n,x}^{(1)}) \leq \sum_{N_x}^{cx^2} \mathbf{P}(V_n > x) \mathbf{P}(\sigma_{-n}^{(1)} > x + n(a - b + \eta_n^{(1)}))$$

$$\leq (1 + o(1)) \frac{d^{(1)}}{a - b} \Delta \overline{F}^s(x).$$

If $n > cx^2$, then
$$\mathbf{P}\left(\frac{V_n}{v\sqrt{n}} > \frac{x}{v\sqrt{n}}\right) = (1 + o(1))2\overline{\Phi}\left(\frac{x}{v\sqrt{n}}\right),$$
as $x \to \infty$, since
$$\mathbf{P}\left(\frac{V_n}{v\sqrt{n}} > y\right) = (1 + o(1))2\overline{\Phi}(y)$$
uniformly in $y$ on a compact set.

Therefore,
$$\sum_{n=cx^2}^{\infty} \mathbf{P}(W^{(2)} > x, A_{n,x}^{(1)})$$
$$= \sum_{cx^2}^{\infty} \mathbf{P}(\sigma_{-n-1}^{(1)} > x + n(a - b + \eta_n))\mathbf{P}(K_{n,x}^{(1,2)})\mathbf{P}(V_n > x, K_{n,x}^{(1,1)})$$
$$= (1 + o(1))d^{(1)} \sum_{cx^2}^{\infty} \overline{F}(x + n(a - b + \eta_n))\left(\mathbf{P}\left(\frac{V_n}{v\sqrt{n}} > \frac{x}{v\sqrt{n}}\right) - o(1)\right)$$
$$= 2(1 + o(1))d^{(1)} \sum_{cx^2}^{\infty} \overline{F}(x + n(a - b))\overline{\Phi}\left(\frac{x}{v\sqrt{n}}\right) + o(\overline{F}^s(x))$$
$$= 2d^{(1)} \int_{cx^2}^{\infty} \overline{F}(x + y(a - b))\overline{\Phi}\left(\frac{x}{v\sqrt{y}}\right) dy + o(\overline{F}^s(x))$$
$$= 2d^{(1)} \int_{0}^{\infty} \overline{F}(x + y(a - b))\overline{\Phi}\left(\frac{x}{v\sqrt{y}}\right) dy$$
$$+ o(\overline{F}^s(x)) - \Delta \frac{d^{(1)}}{a - b}O(\overline{F}^s(x)).$$

Letting $\Delta$ to 0, we get the result. $\square$

*Case $b^{(1)} < b^{(2)}$.*

THEOREM 12. *Assume $b^{(1)} < b^{(2)}$. Then*
$$(57) \quad \mathbf{P}(W^{(2)} > x) = \frac{d^{(2)}}{a - b^{(2)}}\overline{F}^s(x) + \frac{d^{(1)}}{a - b^{(2)}}\overline{F}^s\left(x\frac{a - b^{(1)}}{b^{(2)} - b^{(1)}}\right) + o(\overline{F}^s(x)).$$

*In particular, if either:*

(i) $d^{(2)} > 0$ or
(ii) $d^{(2)} = 0$, $d^{(1)} > 0$ *and*
$$(58) \qquad \liminf_{x \to \infty} \overline{F}^s(2x)/\overline{F}^s(x) > 0,$$

*then one can replace in* (57) *the equality by an equivalence and delete the term* $o(\overline{F}^s(x))$ *in the right-hand side.*

PROOF. Take the notation from the beginning of Section 6.1 and from the proof of Theorem 11. Recall that we consider the case $b^{(2)} = b$. Put $n_x = \frac{x}{b-b^{(1)}}$ and, for a fixed $\varepsilon \in (0,1)$, $n_{x,1} = n_x(1-\varepsilon)$, $n_{x,2} = n_x(1+\varepsilon)$.

Recall that, from (48) and (49), we have to find the asymptotics for

$$\sum_{n \geq N_x} \mathbf{P}(W_0^{(2)} > x, A_{n,x}^{(1)}) = \sum_{n=N_x}^{n_{x,1}} + \sum_{n=n_{x,1}+1}^{n_{x,2}-1} + \sum_{n=n_{x,2}}^{\infty} \equiv P_1(x) + P_2(x) + P_3(x).$$

For $N_x \leq n \leq n_{x,1}$, put $K_{n,x}^{(1)} = \{W_{-n-1}^{(2)} + \sigma_{-n-1}^{(2)} \leq x\}$. Then, for $\eta_n^{(1)} = \eta_n \geq -(a-b)$, on the event $A_{n,x}^{(1)} = K_{n,x}^{(1)} \cap \{\sigma_{-n}^{(1)} > x + n(a-b+\eta_n)\}$,

$$W_{-n}^{(2)} = 0 \quad \text{and} \quad W_0^{(2)} \leq V_n,$$

since $\tau_j \geq \sigma_{j+1}^{(1)}$. Therefore,

$$P_1(x) \leq \sum_{n=N_x}^{n_{x,1}} \mathbf{P}(\sigma_{-n}^{(1)} > x + n(a-b+\eta_n), V_n > x)$$

$$= \sum_{n=N_x}^{n_{x,1}} \mathbf{P}(\sigma_{-n}^{(1)} > x + n(a-b+\eta_n))\mathbf{P}(V_n > x)$$

$$\leq \frac{(1+o(1))d^{(1)}}{a-b}\mathbf{P}(V_{n_{x,1}} > x)\overline{F}^s(x) = o(\overline{F}^s(x)),$$

since $V_n/n \to (b-b^{(1)})^{-1}$ as $n \to \infty$ and

$$\frac{V_{n_{x,1}}}{x} = \frac{V_{n_{x,1}}}{n_{x,1}}\frac{n_{x,1}}{x} \to 1-\varepsilon < 1 \quad \text{a.s.}$$

as $x \to \infty$.

Consider $P_2(x)$. For any sequence $\eta_n \to 0$ and for $x$ sufficiently large,

$$P_2(x) \leq \sum_{n=n_{x,1}+1}^{n_{x,2}-1} \mathbf{P}(\sigma_{-n}^{(1)} > x + n(a-b+\eta_n))$$

$$= \frac{(1+o(1))d^{(1)}}{a-b} \int_{x+n_{x,1}(a-b)}^{x+n_{x,2}(a-b)} \overline{F}(x)\,dx$$

$$\leq \frac{(1+o(1))d^{(1)}}{a-b}\frac{n_{x,2}-n_{x,1}}{n_{x,2}}\overline{F}^s(x) = (1+o(1))\frac{2d^{(1)}\varepsilon}{a-b}\overline{F}^s(x),$$

since $\overline{F}(x)$ is nonincreasing.

Finally, consider $P_3(x)$. We will show that, for the appropriate sequences $\{K_{n,x}^{(1)}\}$ and $\{\eta_n^{(1)}\}$,

(59) $$P_3(x) \sim \frac{d^{(1)}}{a-b}\overline{F}^s\left(x\left(\frac{a-b^{(1)}}{b-b^{(1)}} + \frac{\varepsilon(a-b)}{b-b^{(1)}}\right)\right).$$

Obviously,

$$P_3(x) \leq \sum_{n \geq n_{x,2}} \mathbf{P}(\sigma_{-n}^{(1)} > x + n(a - b + \eta_n^{(1)})),$$

where the right-hand side of the latter inequality is asymptotically equivalent to the right-hand side of (59). Now we establish the lower bound.

From (10), (42), (44) and (45),

$$W_0^{(2)} = Z - W_0^{(1)} - \sigma_0^{(1)} - \sigma_0^{(2)}$$

$$\geq Z_{[-n,0]} - \max_{0 \leq j \leq n} S_j^{(1)} - W_{-n}^{(1)} - \sigma_0^{(1)} - \sigma_0^{(2)}$$

$$\geq \sigma_{-n}^{(1)} - na + \max_{-n \leq q \leq 0}\left(\sum_{m=-n+1}^{q} \sigma_m^{(1)} + S_q^{(\sigma,2)}\right)$$

$$- \max(0, S_{n-1}^{(1)} + \sigma_{-n}^{(1)}) - \max_{0 \leq j \leq n-1} S_j^{(1)} - W_{-n}^{(1)} - \sigma_0^{(1)} - \sigma_0^{(2)}.$$

Due to the SLLN, as $n \to \infty$,

$$\frac{1}{n}S_{n-1}^{(1)} \to b^{(1)} - a, \qquad \frac{1}{n}\left(\max_{0 \leq j \leq n-1} S_j^{(1)} + W_{-n}^{(1)} + \sigma_0^{(1)} + \sigma_0^{(2)}\right) \equiv r_n \to 0$$

and

$$\frac{1}{n}\max_{-n \leq q \leq 0}\left(\sum_{m=-n+1}^{q} \sigma_m^{(1)} + S_q^{(\sigma,2)}\right) \equiv u_n \to b \qquad \text{a.s.}$$

Choose a sequence $\delta_n \downarrow 0$, $n\delta_n \to \infty$ such that, for all $n$,

$$\mathbf{P}(u_n \geq b - \delta_n, r_n \leq \delta_n, S_{n-1}^{(1)} \leq n(b^{(1)} - a + \delta_n)) \geq 1 - \delta_n$$

and denote the latter event by $K_{n,x}^{(1)}$ (it does not depend on $\sigma_{-n}^{(1)}$). On this event,

$$W_0^{(2)} \geq \sigma_{-n}^{(1)} - n(a - b + 2\delta_n) - \max(0, \sigma_{-n}^{(1)} + n(b^{(1)} - a + \delta_n))$$

$$\geq \min(\sigma_{-n}^{(1)} - n(a - b + 2\delta_n), n(b - b^{(1)} - 3\delta_n)).$$

Since $n \geq \frac{(1+\varepsilon)x}{b-b^{(1)}}$, we get

$$n(b - b^{(1)}) - 3n\delta_n \geq x + \frac{\varepsilon(b - b^{(1)})}{1 + \varepsilon}n - 3n\delta_n > x$$

for all sufficiently large $x$. Put $\eta_n^{(1)} = 2\delta_n$. Then, on the event

$$A_{n,x}^{(1)} = \{\sigma_{-n}^{(1)} > x + n(a - b + \eta_n^{(1)})\} \cap K_{n,x}^{(1)},$$

we get $W_0^{(2)} > x$, and (59) follows.

Letting $\varepsilon$ to 0 completes the proof. $\square$

6.2. *Multiserver queues.* The aim of this section is to derive upper and lower bounds and sharp asymptotics for the tail of the stationary maximal dater of multiserver queues. However, we do not obtain here asymptotics for the tail distribution of the stationary waiting time. It is known (see, e.g., [17] and [21]) that these asymptotics may, in general, differ significantly.

Since (AA) does not hold, we cannot use the approach of Section 4. We show how the ideas of Section 5.2 can be used to derive upper and lower bounds which are specific to this queue.

Recall that we can consider a $D/GI/m/\infty$ queue with constant interarrival times $a$. Let $\mathbf{E}\sigma = b$ and $\rho \equiv \frac{b}{ma} \in (0,1)$. Assume further that $\mathbf{P}(\sigma_1 > x) = \overline{F}(x)$, where both distributions $F$ and $F^s$ are subexponential.

THEOREM 13. *Under the foregoing assumptions, when $x$ tends to $\infty$,*

$$\begin{aligned}(60)\quad &\mathbf{P}(Z > x) \\ &= (1 + o(1))\left(\frac{1}{a}\overline{F}^s(x) + \left(\frac{1}{ma - b} - \frac{1}{a}\right)^+ \overline{F}^s\left(\frac{bx}{b - (m-1)a}\right)\right).\end{aligned}$$

Note that the second term in the right-hand side of (60) disappears when $b \leq (m-1)a$.

The proof consists of three steps:

First, we get a lower bound by using the SLLN.
Then we get an upper bound by using results from Section 4.3.1.
Finally, Theorem 7 gives us the tool to derive the exact asymptotics.

LOWER BOUND. Clearly,

$$\mathbf{P}(Z > x) \geq (1 + o(1))\mathbf{P}\left(\bigcup_{n=0}^{\infty} \{\sigma_{-n} > x + na\}\right)$$

$$\sim \sum_n \mathbf{P}(\sigma_{-n} > x + na) \sim \frac{1}{a}\overline{F}^s(x).$$

UPPER BOUND. Take a sufficiently large $L$ and consider the $L$-upperbound $D/GI/1/\infty$ queue with interarrival times $La$ and service times $\{\widehat{s}_n\}$ with mean $\widehat{b} = \mathbf{E}\,\widehat{s}_1$. Since

$$\max_{1 \leq i \leq L} \sigma_i \leq \widehat{s}_1 \leq \sum_{1}^{L} \sigma_i,$$

we get $\mathbf{P}(\widehat{s}_1 > x) \sim L\mathbf{P}(\sigma_1 > x) = L\overline{F}(x)$ as $x \to \infty$. Note that, for the multiserver queue, $\gamma(0) = b/m$. Therefore, we get a natural analogue of Theorem 5,

(61) $$\limsup_{x \to \infty} \frac{\mathbf{P}(Z > x)}{\overline{F}^s(x)} \leq \lim_{L \to \infty} \limsup_{x \to \infty} \frac{\mathbf{P}(\widehat{R} > x)}{\overline{F}^s(x)} = \frac{1}{a - b/m} \overline{F}^s(x).$$

Thus, we are in a position to make use of Theorem 7. The rest of the proof is quite technical and in the same spirit as that of Theorem 9. Because of that, it is omitted.

## APPENDIX

**A.1. Proof of (24).** Put $b = \min(b^{(i_1)}, b^{(i_2)})$. From Corollary 5, we know that, for $r = 1, 2$,

$$\{R^{(i_r)} > x\} \subset \bigcup_{m_r \geq 1} (A_{m_r}^{(i_r)} \cup B^{(i_r)}),$$

where

$$A_{m_r}^{(i_r)} = \{Y_{m_r}^{(i_r)} > x + m_r(b^{(i_r)} - a)\} \quad \text{and} \quad \mathbf{P}(B^{(i_r)}) = o(\overline{F}^s(x)).$$

Therefore,

$$\mathbf{P}(R^{(i_1)} > x, R^{(i_2)} > x)$$

$$\leq \sum_{m_1=1}^{\infty} \sum_{m_2=1}^{\infty} \mathbf{P}(Y_{m_1}^{(i_1)} > x + m_1(b^{(i_1)} - a), Y_{m_2}^{(i_2)} > x + m_2(b^{(i_2)} - a))$$

$$\quad + o(\overline{F}^s(x))$$

(62) $$\leq \sum_{m_1 \neq m_2} \mathbf{P}(Y_1^{(i_1)} > x + m_1(b^{(i_1)} - a)) \mathbf{P}(Y_1^{(i_2)} > x + m_2(b^{(i_2)} - a))$$

$$\quad + \sum_{m=1}^{\infty} \mathbf{P}(\min(Y_1^{(i_1)}, Y_1^{(i_2)}) > x + m(b - a)) + o(\overline{F}^s(x))$$

$$\leq \sum_{m_1} \mathbf{P}(Y_1^{(i_1)} > x + m_1(b^{(i_1)} - a)) \sum_{m_2} \mathbf{P}(Y_1^{(i_2)} > x + m_2(b^{(i_2)} - a))$$

$$\quad + \sum_{m=1}^{\infty} o(\overline{F}(x + mb)) + o(\overline{F}^s(x))$$

$$= \Theta((\overline{F}^s(x))^2) + o(\overline{F}^s(x)) = o(\overline{F}^s(x)).$$

**A.2. Relaxing the independence assumptions.** The aim of this section is to give conditions under which assumption (H) of Section 4.1 is satisfied, although the r.v.'s $Y^{(j)}$ are not independent.

We assume that there exists a random variable $\nu$ taking values in an arbitrary measurable space $(\mathcal{Y}, B_Y)$ and such that:

- Given $\nu$, the random variables $Y_1^{(j)}, j = 1, \ldots, r$, are conditionally independent.
- For any $j = 1, \ldots, r$,

$$\mathbf{P}(Y_1^{(j)} > x|\nu) \sim d_\nu^{(j)} \overline{F}(x), \tag{63}$$

$\mathbf{P}_\nu$-a.s., where $d_\nu^{(j)}$ is a nonnegative random variable with a finite mean $d^{(j)}$.

Then

$$\widetilde{d}_\nu^{(j)} \equiv \sup_x \frac{\mathbf{P}(Y_1^{(j)} > x|\nu)}{\overline{F}(x)} \tag{64}$$

is an a.s. finite random variable, too.

Assume in addition that, for any $1 \leq j_1 \leq j_2 \leq r$,

$$\mathbf{E} \prod_{j=j_1}^{j_2} \widetilde{d}_\nu^{(j)} < \infty. \tag{65}$$

LEMMA 8. *Under the foregoing assumptions, for any $1 \leq j_1 \leq j_2 \leq r$,*

$$\mathbf{P}\bigg(\sum_{j=j_1}^{j_2} Y_1^{(j)} > x\bigg) \sim \mathbf{P}\bigg(\max_{j_1 \leq j \leq j_2} Y_1^{(j)} > x\bigg)$$
$$\sim \sum_{j=j_1}^{j_2} \mathbf{P}(Y_1^{(j)} > x) \sim \sum_{j=j_1}^{j_2} d^{(j)} \overline{F}(x). \tag{66}$$

PROOF. Without loss of generality, we prove the result for $j_1 = 1$, $j_2 = r$. Note that

$$\frac{\mathbf{P}(\sum_1^r Y_1^{(j)} > x|\nu)}{\overline{F}(x)} \to \sum_1^r d_\nu^{(j)} \leftarrow \frac{\mathbf{P}(\max_j Y_1^{(j)} > x|\nu)}{\overline{F}(x)},$$

$\mathbf{P}_\nu$-a.s. and, for all $x$,

$$0 \leq \frac{\mathbf{P}(\max_j Y_1^{(j)} > x|\nu)}{\overline{F}(x)} \leq \frac{\mathbf{P}(\sum_1^r Y_1^{(j)} > x|\nu)}{\overline{F}(x)} \leq \prod_1^r \widetilde{d}_\nu^{(j)} \cdot \sup_x \frac{\overline{F}^{*r}(x)}{\overline{F}(x)},$$

where the latter supremum is finite. Then the dominated convergence theorem implies that

$$\frac{\mathbf{P}(\sum_1^r Y^{(j)} > x)}{\overline{F}(x)} \equiv \mathbf{E}\left(\frac{\mathbf{P}(\sum_1^r Y_1^{(j)} > x|\nu)}{\overline{F}(x)}\right) \to d \equiv \sum_j \mathbf{E} d_\nu^{(j)} \equiv \sum_j d^{(j)}$$

and

$$\frac{\mathbf{P}(\max_j Y_1^{(j)} > x)}{\overline{F}(x)} \to d.$$

□

Consider the following example, which covers the generalized Jackson network case. Assume that there are given:

- Some random vector $\nu = (\nu^{(1)}, \ldots, \nu^{(r)})$ with nonnegative integer-valued components, such that $\mathbf{E}\exp(c\nu^{(j)}) < \infty$ for some $c > 0$ and for all $j = 1, \ldots, r$;
- $r$ sequences $\{\sigma_n^{(j)}\}$ of i.i.d. subexponential random variables that are mutually independent and independent of $\nu$, and such that $\mathbf{P}(\sigma_1^{(j)} > x) \sim l^{(j)}\overline{F}(x)$. We do not make the assumption that the r.v.'s $\nu^{(1)}, \ldots, \nu^{(r)}$ are independent.

Put $Y_1^{(j)} = \sum_{i=1}^{\nu^{(j)}} \sigma_i^{(j)}$. The above conditions imply that

$$\mathbf{E}\exp\left(\frac{c}{r}\sum_j \nu^{(j)}\right) \leq \mathbf{E}\exp\left(c\max_j \nu^{(j)}\right) \leq \sum_j \mathbf{E}\exp(c\nu^{(j)}) < \infty$$

and that for all $j = 1, \ldots, r$,

$$u^{(j)} \equiv \sup_t \frac{\mathbf{P}(\sigma_1^{(j)} > t)}{\overline{F}(t)} < \infty.$$

Due to subexponentiality, for $j = 1, \ldots, r$,

$$\mathbf{P}(Y_1^{(j)} > x|\nu) \sim \nu^{(j)} l^{(j)} \overline{F}(x).$$

It is known (see, e.g., [13], page 41) that, for any $\varepsilon > 0$, one can choose $K^{(j)} \equiv K^{(j)}(\varepsilon)$ such that

$$\mathbf{P}(Y_1^{(j)} > x|\nu) \leq K^{(j)}(1+\varepsilon)^{\nu^{(j)}} \mathbf{P}(\sigma_1^{(j)} > x).$$

The right-hand side of the latter inequality is not bigger than $K^{(j)} u^{(j)}(1+\varepsilon)^{\nu^{(j)}} \overline{F}(x)$.

Take $\varepsilon > 0$ such that $\log(1+\varepsilon) \leq c'$. Then the conditions of Lemma 8 are satisfied with $d_\nu^{(j)} = \nu^{(j)} l^{(j)}$, $d^{(j)} = l^{(j)} \mathbf{E}\nu^{(j)}$ and $\widetilde{d}_\nu^{(j)} = K^{(j)} u^{(j)} (1+\varepsilon)^{\nu^{(j)}}$.

**A.3. Deterministic interarrival times.** We extend to the monotone separable framework the approach used in [4] for single server queues to show that there may be no loss of generality in assuming that a network has deterministic interarrival times when one wants to evaluate the tail asymptotics of its stationary maximal dater.

The framework is that of Section 2. Fix $\{\zeta_n\}$, $\{f_l\}$ and consider a family of networks with different "input sequences" $\{T_n\}$ such that $\mathbf{E}\tau_1 > \gamma(0)$. Without loss of generality assume $T_0 = 0$.

In particular, a network with constant interarrival times (say $a$) belongs to this family. For such a network, we use the notation $Z^{(a)}$ and $Z^{(a)}_{[-n,0]}$.

For any $\{T_n\}$ and for any $\widetilde{a} < \mathbf{E}\tau_1$, set

$$\psi(\{\tau_n\}, \widetilde{a}) = \sup_{n \geq 0}(n\widetilde{a} + T_{-n}) \equiv \sup_{n \geq 0}\left(\sum_{i=-n}^{-1}(\widetilde{a} - \tau_i)\right).$$

THEOREM 14. *Assume that there exist a continuous and strictly positive function $h:(\gamma(0), \infty) \to (0, \infty)$ and a subexponential distribution $G$ such that, for any $a > \gamma(0)$,*

$$(67) \qquad \mathbf{P}(Z^{(a)} > x) \sim h(a)\overline{G}(x) \qquad \text{as } x \to \infty.$$

*Then, for any network with random interarrival times $\{\tau_n\}$, such that $\mathbf{E}\tau_1 = a > \gamma(0)$, the following is valid: If $\{\tau_n\}$ and $\{\zeta_n\}$ are independent and if, for any $\widetilde{a} < a$,*

$$(68) \qquad \mathbf{P}(\psi(\{\tau_n\}, \widetilde{a}) > x) = o(\overline{G}(x)) \qquad \text{as } x \to \infty,$$

*then*

$$(69) \qquad \mathbf{P}(Z > x) \sim h(a)\overline{G}(x) \qquad \text{as } x \to \infty.$$

REMARK 8. In particular, condition (68) is satisfied if the $\tau_n$'s are i.i.d. Indeed, then $\psi(\{\tau_n\}, \widetilde{a})$ has either a bounded [if $\mathbf{P}(\widetilde{a} \geq \tau_1) = 0$] or exponential tail, which is lighter than any long tail.

PROOF OF THEOREM 14. Take any $\varepsilon \in (0, a - \gamma(0))$. Due to the monotonicity,

$$Z_{[-n,0]} \leq Z^{(a-\varepsilon)}_{[-n,0]} + \max_{0 \leq j \leq n}\left(\sum_{i=-j}^{-1}(a - \varepsilon - \tau_i)\right)^+.$$

Therefore,

$$Z \leq Z^{(a-\varepsilon)} + \psi(\{\tau_n\}, a - \varepsilon) \equiv Z^{(a-\varepsilon)} + \psi,$$

where $Z^{(a-\varepsilon)}$ and $\psi$ are independent. Therefore,

$$\mathbf{P}(Z > x) \leq \mathbf{P}(Z^{(a-\varepsilon)} + \psi > x) \sim \mathbf{P}(Z^{(a-\varepsilon)} > x) \sim h(a-\varepsilon)\overline{G}(x).$$

Thus

$$\limsup_{x \to \infty} \frac{\mathbf{P}(Z > x)}{\overline{G}(x)} \leq h(a-\varepsilon)$$

for any $\varepsilon \in (0, a - \gamma(0))$. Letting $\varepsilon$ go to 0, we get the upper bound $h(a)$.

For the lower bound, we use the monotonicity, the SLLN for the $\tau$'s, the LT and the independence assumptions. For any $\varepsilon > 0$, one can choose a sufficiently large $C \equiv C(\varepsilon)$ such that

$$\mathbf{P}(T_{-n} \geq -n(a+\varepsilon) - C \; \forall n \geq 0) \geq 1 - \varepsilon.$$

Denote the latter event by $D_\varepsilon$. Then

$$\begin{aligned}\mathbf{P}(Z > x) &\geq \mathbf{P}(Z > x, D_\varepsilon) \geq \mathbf{P}(Z^{(a+\varepsilon)} - C > x, D_\varepsilon) \\ &\geq \mathbf{P}(Z^{(a+\varepsilon)} - C > x)(1-\varepsilon) \sim h(a+\varepsilon)(1-\varepsilon)\overline{G}(x+C) \\ &\sim h(a+\varepsilon)(1-\varepsilon)\overline{G}(x).\end{aligned}$$

Thus, for any $\varepsilon \in (0,1)$,

$$\liminf_{x \to \infty} \frac{\mathbf{P}(Z > x)}{\overline{G}(x)} \geq h(a+\varepsilon)(1-\varepsilon).$$

Letting $\varepsilon$ go to 0, we get the lower bound with coincides with the upper one. □

INRIA-ENS  
45 RUE D'ULM  
75005 PARIS  
FRANCE  
E-MAIL: Francois.Baccelli@ens.fr

DEPARTMENT OF ACTUARIAL MATHEMATICS  
AND STATISTICS  
HERIOT-WATT UNIVERSITY  
EDINBURGH  
UNITED KINGDOM  
E-MAIL: S.Foss@ma.hw.ac.uk